\documentclass[12pt]{article}
\usepackage{amsmath,amsfonts,amsthm,mathrsfs}

\newcommand{\R}{{\mathbb R}} \newcommand{\N}{{\mathbb N}}
\newcommand{\K}{{\mathbb K}} \newcommand{\Z}{{\mathbb Z}}
  \def\C{{\mathbb C}}
\newcommand{\Prm}{{\mathbb P}}
\newcommand{\wt}{\widetilde }
\renewcommand{\epsilon}{\varepsilon } 
 
\renewcommand{\rho}{\varrho } 
\renewcommand{\phi}{\varphi }
\newcommand{\E}{{\mathbb E}\,}
\newcommand{\EE}{{\mathbb E}}

\newcommand{\ran}{{\rm ran }}
\newcommand{\set}{{\rm set }}

\newcommand{\ranno}{\text{\rm ran--non}}
\newcommand{\de}{{\rm det }}
\newcommand{\deno}{\text{\rm det--non}}
\newcommand{\ca}{{\rm card}}
\newcommand{\avg}{{\rm avg }}
\newcommand{\avgno}{\text{\rm avg--non}}

\newcommand{\SUPP}{{\rm supp}}

\newtheorem{theorem}{Theorem}[section]
\newtheorem{lemma}[theorem]{Lemma}
\newtheorem{corollary}[theorem]{Corollary}
\newtheorem{proposition}[theorem]{Proposition}

\newtheorem{remark}[theorem]{Remark}

\begin {document}
 \title{
Randomized Complexity of Mean Computation and the Adaption Problem
}

 \author {Stefan Heinrich\\
Department of Computer Science\\
RPTU Kaiserslautern-Landau\\
D-67653 Kaiserslautern, Germany}  
\date{\today}
\maketitle

\begin{abstract} 
Recently the adaption problem of Information-Based Complexity (IBC) for linear problems in the randomized setting was solved in Heinrich (2024)\cite{Hei23a}. Several papers treating further aspects of this problem followed. However, all examples obtained so far were vector-valued. In this paper we settle the scalar-valued case. We study the complexity of mean computation in finite dimensional sequence spaces with mixed $L_p^N$ norms.
We determine the $n$-th minimal errors in the randomized adaptive and non-adaptive setting.  It turns out that among the problems considered there are examples where adaptive and non-adaptive $n$-th minimal errors deviate by a power of $n$. The gap can be (up to log factors) of the order $n^{1/4}$. We also show how to turn such results into infinite dimensional examples with suitable deviation for all $n$ simultaneously.
\end{abstract}

\section{Introduction}
\label{sec:1}
Let $N,N_1,N_2\in{\mathbb{N}}$,  $1\leq p,u \leq \infty$, and  let
$L_p^N$ be the space of all functions $f: {\mathbb{Z}}[1,N]:=\{1,2,\dots,N\} \rightarrow {\mathbb{K}}$
with the norm
\begin{eqnarray*}
\| f \|_{L_p^N} = 
  \left( \frac{1}{N} \sum_{i=1}^N |f(i)|^p \right)^{1/p} \;(p<\infty),\quad  
  \| f \|_{L_\infty^N}=\max_{1\le i\le N} |f(i)| .
\end{eqnarray*}
Define the space $L_p^{N_1}\big(L_u^{N_2}\big)$ as the set of all functions $f: {\mathbb{Z}}[1,N_1]\times {\mathbb{Z}}[1,N_2]\rightarrow {\mathbb{K}}$
endowed with the norm
\[
\| f \|_{L_p^{N_1}\big(L_u^{N_2}\big)} = \Big\|\big(\|f_i\|_{L_u^{N_2}}\big)_{i=1}^{N_1}\Big\|_{L_p^{N_1}},
\]
where 
$
f_i=(f(i,j))_{j=1}^{N_2}
$
are the rows of the matrix $(f(i,j))$.
In the present paper we study the complexity of mean computation in the randomized setting.
We determine the order of the randomized $n$-th minimal errors of 
\begin{equation}
\label{C6}
I^{N_1,N_2}:L_p^{N_1}\big(L_u^{N_2}\big)\to \K, \quad I^{N_1,N_2}f=\frac{1}{N_1N_2} \sum_{i=1}^{N_1} \sum_{j=1}^{N_2}f(i,j).
\end{equation}
The input set is the unit ball of $L_p^{N_1}\big(L_u^{N_2}\big)$ and information is standard (values of $f$).

The adaption problem of Information-Based Complexity (IBC) for linear problems is concerned with the relation between adaptive and non-adaptive $n$-th minimal errors. In 1980 Gal, Micchelli \cite{GM80} and  Traub, Wo\'zniakowski \cite{TW80} showed that in the deterministic setting  adaptive and non-adaptive deterministic $n$-th minimal errors can deviate at most by a factor of 2: for any linear problem $\mathcal{P}=(F,G,S,K,\Lambda)$ and any $n\in {\mathbb{N}}$
\begin{equation}
\label{J9}
e_n^{\rm det-non } (S,F,G)\le 2e_n^{\rm det } (S,F,G).
\end{equation}
In 1996 Novak \cite{Nov96} posed the respective problem for the randomized setting:
Is there a constant $c>0$ such that for all linear problems
$\mathcal{P}=(F,G,S,K,\Lambda)$ and all $n\in{\mathbb{N}}$
\begin{equation}
e_n^{\rm ran-non} (S,F,G)\le ce_n^{\rm ran } (S,F,G) \, {\bf ?}
\label{A0}
\end{equation}
See the problem on p.\ 213 of \cite{Nov96}, and also Problem 20 on p.\ 146 of the monograph \cite{NW08} by Novak and Wozniakowski (2008).
This problem was solved recently by the author in \cite{Hei23a}, where it was shown that for some instances of vector-valued mean computation the gap between non-adaptive and adaptive randomized $n$-th minimal errors can be (up to log factors) of  order $n^{1/8}$.  Considering vector valued approximation, it was shown in another paper by author \cite{Hei23c} that the gap can be $n^{1/2}$ (again, up to log factors). Both papers deal with standard information, that is, function values. Problem \eqref{A0} remained open for the case of arbitrary linear information (that is, $\Lambda$ consists of all linear functionals on $F$). This was settled recently by Kunsch, Novak, and Wnuk \cite{KNW23}. 

All counter-examples to  problem \eqref{A0} given so far were vector valued. The scalar-valued case $G=\K$ remained open. In the present paper we show that the answer is negative, as well. 
In the case $1\le p<2<u\le \infty$ of mean computation \eqref{C6} adaptive and non-adaptive randomized $n$-th minimal errors deviate by a power of $n$, see 
relations \eqref{A3} and \eqref{A4} of Theorem \ref{theo:1}. This is done by showing that for each $n$ there is a finite dimensional integration problem so that the gaps increase with growing $n$.

This raises the question about infinite dimensional examples with respective gaps for all $n$ simultaneously. For vector-valued mean computation such an example - namely parametric integration - was presented in \cite{Hei23b}. We show that there are such infinite dimensional examples also for integration. For this purpose we use an  approach different from that in \cite{Hei23b}. We present a general way of passing from finite into infinite dimensional examples by the help of direct sums.

The paper is organized as follows. In Section \ref{sec:2} we recall the basic notions of IBC and present some auxiliary facts. Moreover, this section contains a new general result on the average case setting for sum problems. Section \ref{sec:3} presents non-adaptive and adaptive algorithms for mean computation and their error estimates. Lower bounds and the main complexity result are contained in Section \ref{sec:4}, while Section \ref{sec:5} is devoted to the procedure of passage to infinite dimensional problems for mean computation (Subsection \ref{sec:5.1}) as well as,  based on the results of \cite{Hei23c}, for approximation (Subsection \ref{sec:5.2}).

\section{Preliminaries}
\label{sec:2}

We denote $\N=\{1,2,\dots\}$, $\N_0=\N\cup\{0\}$, and $\Z[N_1,N_2]=\{1,2,\dots ,N\}$ for $N_1,N_2\in\N_0$, $N_1\le N_2$. The symbol $\K$ stands for the scalar field $\R$ or $\C$.
We often use the same symbol
$c, c_1,c_2,\dots$ for possibly different constants, even if they appear in a sequence
of relations. However, some constants are supposed to have the same meaning throughout a proof  --  these are denoted by symbols $c(1),c(2),\dots$. The unit ball of a normed space $X$ is denoted by $B_X$. Throughout the paper $\log$ means $\log_2$. 

We adopt the general IBC notation as presented in in Section 2 of \cite{Hei23a}. For background and all details we refer to \cite{Hei23a} as well as to  \cite{Nov88,TWW88} and \cite{Hei05a, Hei05b}. 
An abstract  numerical problem $\mathcal{P}$ is given as 
\begin{equation}
\label{M7}
\mathcal{P}=(F,G,S,K,\Lambda).
\end{equation}
where $F$ is a non-empty set, 
$G$ a Banach space and $S$ is a mapping $F\to G$. Furthermore, $K$ is any nonempty set and $\Lambda$ is a nonempty set of mappings from $F$ to $K$. The operator $S$ is called the  solution operator and $\Lambda$ the set of information functionals. 
A problem $\mathcal{P}$ is called linear, if $K=\K$, $F$ is a convex and balanced subset of a linear space $X$ over $\K$,  
$S$ is the restriction to $F$ of a linear operator 
from $X$ to $G$, and each $\lambda\in\Lambda$ is the restriction  to $F$ of a linear mapping from $X$ to $\K$. 

For $n\in\N_0$ the adaptive (respectively non-adaptive) deterministic $n$-th minimal error of $S$ is denoted by $e_n^\de(S,F,G)$  ($e_n^\deno(S,F,G)$). Correspondingly, 
$e_n^\ran(S,F,G)$  ($e_n^\ranno(S,F,G)$  stand for the adaptive (respectively non-adaptive) randomized $n$-th minimal error of $S$. Furthermore, given a probability measure $\mu$ on $F$ whose support is a finite set, $e_n^{\rm avg }(S,\mu,G)$  ($e_n^\avgno(S,\mu,G))$ denote the adaptive (respectively non-adaptive)  $n$-th minimal average error of $S$.

The following relations hold for $n\in \N_0$
\begin{eqnarray}
e_n^\de (S,F,G)&\le& e_n^\deno (S,F,G)
\label{AF5}\\
e_n^\ran (S,F,G)&\le& e_n^\ranno (S,F,G)
\label{AF6}\\
e_n^\ran (S,F,G)&\le&e_n^\de (S,F,G)
\label{G5}\\
e_n^\ranno (S,F,G)&\le& e_n^\deno (S,F,G)
\label{AF7}\\
e_n^\avg (S,\mu,G)&\le& e_n^\avgno (S,\mu,G)\label{B7}
\end{eqnarray}
and for each probability measure $\mu$ on $F$ of finite support 
\begin{eqnarray}
e_n^\ran(S,F,G)&\ge& \frac{1}{2}e_{2n}^\avg(S,\mu,G)\label{CK3}
\\
e_n^\ranno(S,F,G)&\ge& \frac{1}{2}e_{2n}^\avgno(S,\mu,G).\label{CK4}
\end{eqnarray}

Similar to \cite{Hei23a}, Section 2,  we need some further general results on algorithms in product structures. Let $M\in \N$ and let
$\mathcal{P}_i=(F_i,G,S_i,K_i,\Lambda_i)$ $(i=1,\dots, M)$ be numerical problems (with the same target space $G$ for all $i$). We assume that for each $i$ none of the elements of $\Lambda_i$ is constant on $F_i$, that is, 
\begin{equation}
\label{J1}
\text{for all\;} \lambda\in\Lambda_i\; \text{there exist\;}f_1,f_2\in F_i\;\text{with\;}\lambda(f_1)\ne \lambda(f_2).
\end{equation}
Define the sum problem $\mathcal{P}=(F,G,S,K,\Lambda)$ by 
\begin{eqnarray*}
F=\prod_{i=1}^M F_i,
\quad S:F\to G,\; S(f_1,\dots,f_M)=\sum_{i=1}^M S_i(f_i),\quad K=\bigcup_{i=1}^M K_i,\quad \Lambda=\bigcup_{i=1}^M\Phi_i(\Lambda_i), 
\end{eqnarray*}
where
\begin{equation*}
\Phi_i:\Lambda_i\to \mathscr{F}(F,K), \quad (\Phi_i(\lambda_i))(f_1,\dots,f_i,\dots, f_M)=\lambda_i(f_i),
\end{equation*}
and $\mathscr{F}(F,K)$ stands for the set of all mappings from $F$ to $K$.
Observe that \eqref{J1} implies
$\Phi_i(\Lambda_i)\cap\Phi_j(\Lambda_j)=\emptyset$ $(i\ne j)$.
For $1\le i\le M$ we set
\begin{equation*}
F'_i=\prod_{1\le j\le M,j\ne i} F_j.
\end{equation*}
If $i$ is fixed, we  identify, for convenience of notation, 
\begin{equation*}
F\quad \text{with}\quad F_i\times F_i', \quad f=(f_1,\dots,f_i,\dots,f_M)\in F \quad \text{with}\quad f=(f_i,f'_i),
\end{equation*}
 where 
\begin{equation*}
f'_i=(f_1,\dots,f_{i-1},f_{i+1},\dots, f_M)\in F_i'.
\end{equation*}
Let $f_{i,0}'=(f_{i,1,0},\dots,f_{i,i-1,0},f_{i,i+1,0},\dots, f_{i,M,0})\in F_i'$ be fixed elements with the property
\begin{equation}
\label{P1}
\sum_{1\le j\le M,j\ne i}S_j(f_{i,j,0})=0,
\end{equation}
and let
\begin{equation}
J_i:F_i\to F,\quad J_i(f_i)=(f_i,f_{i,0}')\quad (f_i\in F_i).
\label{P6}
\end{equation}
Now let $\mu_i$ be probability measures on $F_i$ whose support is a finite set and let $\nu_i\ge 0$ be reals with $\sum_{i=1}^M \nu_i=1$. 
We define a measure $\mu$ on $F$ of finite support by setting for a set $C\subset F$
\begin{equation}
\label{B8}
\mu(C)=\sum_{i=1}^M\nu_i\mu_i(J_i^{-1}(C)).
\end{equation}
\begin{lemma}
\label{Ulem:3}
With the notation above and under assumption \eqref{J1} we have for each $n\in\N_0$
\begin{equation}
\label{UN5}
e_n^\avgno(S,\mu,G)\ge \min\Bigg\{\sum_{i=1}^M\nu_ie_{n_i}^\avgno(S_i,\mu_i,G):\,n_i\in\N_0, n_i\ge 0,\sum_{i=1}^Mn_i\le n\Bigg\}.
\end{equation}
\end{lemma}
\begin{proof}
Let $A$ be a non-adaptive deterministic algorithm for $\mathcal{P}$ with $\ca(A)\le n$.  Let $n_i$ be the number of those information functionals of $A$ which are from $\Phi_i(\Lambda_i)$. Then 
\begin{equation}
\label{UN8}
\sum_{i=1}^Mn_i\le n.
\end{equation}
For fixed $i$ we will use Lemma 2.1 of \cite{Hei23a} with
\begin{eqnarray*}
&&F^{(1)}=F_i,\quad F^{(2)}=F_i',\quad K^{(1)}=K_i,\quad  K^{(2)}=\bigcup_{j:j\ne i} K_j,
\\
&&\Lambda^{(1)}=\Phi_i(\Lambda_i), \quad\Lambda^{(2)}=\bigcup_{j:j\ne i} \Phi_j(\Lambda_j).
\end{eqnarray*}
Let 
\begin{equation*}
\mathcal{P}_{f_{i,0}'}=(F_i,G,S_{f_{i,0}'},K_i,\Lambda_{f_{i,0}'}),
\end{equation*}
be the restricted problem obtained by fixing the second component to be $f_{i,0}'$, that is,
\begin{eqnarray*}
S_{f_{i,0}'}(f_i)=S(f_i,f_{i,0}')\quad(f_i\in F_i),
\quad
\Lambda_{f_{i,0}'}=\{\lambda(\,\cdot\,,f_{i,0}'):\; \lambda\in \Phi_i(\Lambda_i)\}.
\end{eqnarray*}
Now Lemma 2.1 of \cite{Hei23a} shows that there is a deterministic non-adaptive algorithm $A_{i,f'_{i,0}}$ for $\mathcal{P}_{f_{i,0}'}$ such that for all $f_i\in F_i$ 
\begin{eqnarray}
A_{i,f_{i,0}'}(f_i)&=&A(f_i,f_{i,0}')
\label{N4}\\
\ca(A_{i,f_{i,0}'})&=&n_i.
\label{N5}
\end{eqnarray}
Moreover, by \eqref{P1}, for $f_i\in F_i$ 
\begin{equation}
S_{f_{i,0}'}(f_i)=S(f_i,f_{i,0}')=S_i(f_i)+\sum_{j: j\ne i}S_j(f_{i,j,0})=S_i(f_i),
\label{P5}
\end{equation}
and, since for  $\lambda_i\in\Lambda_i$ we have $(\Phi_i(\lambda_i))(f_i,f_{i,0}')=\lambda_i(f_i)$,
\begin{equation*}
\Lambda_{f_{i,0}'}=\{\lambda(\,\cdot\,,f_{i,0}'):\; \lambda\in \Phi_i(\Lambda_i)\}
=\{(\Phi_i(\lambda_i))(\,\cdot\,,f_{i,0}'):\; \lambda_i\in \Lambda_i\}=\Lambda_i.
\end{equation*}
This implies
$\mathcal{P}_{f_{i,0}'}=\mathcal{P}_i$,
so $A_{i,f_{i,0}'}$ is a deterministic non-adaptive algorithm for $\mathcal{P}_i$. 
From \eqref{P6}, \eqref{N4}, and \eqref{P5} we conclude
\begin{eqnarray*}
A(J_i(f_i))=A(f_i,f_{i,0}')=A_{i,f_{i,0}'}(f_i),\quad   
S(J_i(f_i))=S(f_i,f_{i,0}')=S_i(f_i).
\end{eqnarray*}
Consequently, using also \eqref{B8} and \eqref{N5},
\begin{eqnarray*}
\int_F\|S(f)-A(f)\|_G d\mu(f)&=&\sum_{i=1}^M \nu_i\int_{F_i}\|S(J_i(f_i))-A(J_i(f_i))\|_Gd\mu_i(f_i)
\nonumber\\
&=&\sum_{i=1}^M \nu_i\int_{F_i}\|S_i(f_i)-A_{i,f'_{i,0}}(f_i)\|_G d\mu_i(f_i)
\nonumber\\
&\ge &\sum_{i=1}^M \nu_ie_{n_i}^\avgno(S_i,\mu_i,G),
\end{eqnarray*}
which together with \eqref{UN8} implies \eqref{UN5}.

\end{proof}
Now consider the case that all $\mathcal{P}_i$ are copies of the same problem $\mathcal{P}_1=(F_1,G,S_1,K_1,\Lambda_1)$,  $\mu_i=\mu_1$, $\nu_i=M^{-1}$ $(i=1,\dots,M)$.
\begin{corollary}
\label{Ucor:2}
\begin{equation}
\label{U0}
e_n^\avgno(S,\mu,G)\ge 2^{-1}e_{\left\lfloor\frac{2n}{M}\right\rfloor}^\avgno(S_1,\mu_1,G).
\end{equation}
\end{corollary}
\begin{proof}
Let $n_i\in\N_0$,   $\sum_{i=1}^M n_i\le n$ and define $\mathcal{I}=\{i:\, n_i\le \frac{2n}{M}\}$, thus $|\mathcal{I}|\ge \frac{M}{2}$. Hence, for $i\in \mathcal{I}$, 
\begin{equation*}
e_{n_i}^\avgno(S_1,\mu_1,G)\ge e_{\left\lfloor\frac{2n}{M}\right\rfloor}^\avgno(S_1,\mu_1,G),
\end{equation*}
so the desired result follows from Lemma \ref{Ulem:3}.

\end{proof}
Next we state a result on reduction.
Let $\mathcal{P}=(F,G,S,K,\Lambda)$ and $\wt{\mathcal{P}}=(\wt{F},\wt{G},\wt{S},\wt{K},\wt{\Lambda})$ be numerical problems. We say that 
$\mathcal{P}$ reduces to $\wt{\mathcal{P}}$, if the following holds. There are mappings $R:F\to\wt{F}$ and $\Psi:\wt{G}\to G$ such that 
$$
S=\Psi\circ\wt{S}\circ R. 
$$
Furthermore, there exist a $\kappa \in \N$, mappings 
$\eta_j:\wt{\Lambda}\to\Lambda \quad (j=1, \dots, \kappa)$
and $\rho:\wt{\Lambda}\times K^\kappa\to \wt{K}$ such that
\begin{equation}
\label{XG3}
\wt{\lambda}(R(f))=\rho\Big(\wt{\lambda}, \big(\eta_1(\wt{\lambda})\big)(f), \dots,\big(\eta_\kappa(\wt{\lambda})\big)(f)\Big)
\end{equation}
for all $f\in F$ and $\wt{\lambda}\in\wt{\Lambda}$. Finally, we assume that $\Psi:\wt{G}\to G$ is a Lipschitz mapping, that is, there is a constant $c\ge 0$ such that
$$
\|\Psi(x)-\Psi(y)\|_G \le c\,\|x-y\|_{\wt{G}}\quad\mbox{for all}\quad x,y\in\wt{G}.
$$ 
The Lipschitz constant $\|\Psi\|_{\rm Lip}$ is the smallest constant $c$ such that
the relation above holds.  We refer to \cite{Hei05b}, Section 3 for this notion and some background.
The following is Proposition 1 of \cite{Hei05b}.
\begin{proposition}
\label{pro:2}  
Suppose that $\mathcal{P}$ reduces to $\wt{\mathcal{P}}$ and let ${\rm set}\in\{\de,\ran,\ranno,\deno\}$. Then for all $n\in \N_0$,
\begin{eqnarray*}
e_{\kappa n}^\set(S,F,G)&\le& \|\Psi\|_{\rm Lip}\, e_n^\set(\wt{S},\wt{F},\wt{G}).
\end{eqnarray*}
\end{proposition}
The following result, which is Proposition 2 of \cite{Hei05b}, concerns 
additivity properties of the minimal errors. 
\begin{proposition}
\label{pro:3}
 Let  ${\rm set}\in\{\de,\deno,\ran,\ranno\}$, $M\in\N$
and let $S_k:F\to G$ $(k=1,\dots,M)$ be mappings. Define $S:F\to G$
by $S(f)=\sum_{k=1}^M S_k(f)\quad(f\in F)$. 
 Let $n_1,\dots,n_M\in\N_0$ and put
$n=\sum_{k=1}^M n_k$. 
Then
\begin{eqnarray*}
e_n^\set(S,F)&\le&\sum_{k=1}^M e_{n_k}^\set(S_k,F).
\end{eqnarray*}
\end{proposition}
Propositions \ref{pro:2} and \ref{pro:3}  were proved in \cite{Hei05b} for the adaptive setting. This is the technically involved case. The non-adaptive case is much easier, essentially straight-forward. We omit the proofs.
The  next lemma is well-known in IBC, see \cite{Nov88,TWW88}, and  
specifically  \cite{Hei05a},  Lemma 6 for statement (i) and \cite{Hei18a}, Proposition 3.1 for (ii). 
\begin{lemma}\label{lem:5} 
Assume that  $K=\K$, $F$ is a  subset of a linear space $X$ over $\K$,  
$S$ is the restriction to $F$ of a linear operator 
from $X$ to $G$, and each $\lambda\in\Lambda$ is the restriction to $F$  of a linear mapping from $X$ to $\K$.
Let $\bar{n}\in\N$ and suppose there are $(\psi_i)_{i=1}^{\bar{n}}\subseteq F$
such that the sets $\{\lambda\in \Lambda\,:\, \lambda(\psi_i)\ne 0\}\; (i=1,\dots,\bar{n})$
are mutually disjoint.
Then the following hold for all $n\in\N_0$ with
$4n<\bar{n}$:

(i) If $\sum_{i=1}^{\bar{n}} \alpha_i \psi_i\in F$ for all sequences $(\alpha_i)_{i=1}^{\bar{n}}\in \{-1,1\}^{\bar{n}}$ and $\mu$ is the distribution of $\sum_{i=1}^{\bar{n}} \epsilon_i \psi_i$, where $\epsilon_i$ are independent Bernoulli random variables
with $\Prm\{\epsilon_i=1\}=\Prm\{\epsilon_i=-1\}=1/2$,  then 
$$
e_n^\avg(S,\mu,G)\ge \frac{1}{2}\min\bigg\{\E\Big\|\sum_{i\in \mathcal{I}}\epsilon_iS\psi_i\Big\|_G:\,\mathcal{I}\subseteq\{1,\dots,\bar{n}\},\,|\mathcal{I}|\ge \bar{n}-2n\bigg\}.
$$

(ii) If $\alpha \psi_i\in F$ for all $1\le i\le\bar{n}$ and $\alpha\in \{-1,1\}$, and $\mu$ is the uniform distribution on the set $\{\alpha \psi_i\,:\,1\le i\le\bar{n},\; \alpha\in \{-1,1\}\}$,
then 
$$
e_n^\avg(S,\mu,G)\ge \frac{1}{2}\min_{1\le i\le \bar{n}} \|S\psi_i\|_G.
$$
\end{lemma}

Let $\theta$ be the mapping given by the median, that is, if $z_1^*\le  \dots\le z_m^*$ is the non-decreasing rearrangement of $(z_1,\dots,z_m)\in{\mathbb R}^m$, then  
$\theta(z_1,\dots,z_m)$ stands for 
$z^*_{(m+1)/2}$ if $m$ is odd and $\frac{z_{m/2}^*+z_{m/2+1}^*}{2}$ if $m$ is even. The following is well-known,  see, e.g, \cite{{Hei01}}.
\begin{lemma}
\label{Ulem:2e}
Let $\zeta_1,\dots,\zeta_m$ be independent, identically distributed real-valued random variables on a probability space $(\Omega,\Sigma,{\mathbb P})$, $z\in {\mathbb R}$, $\varepsilon>0$ , and assume that ${\mathbb P}\{|z-\zeta_1|\le\varepsilon\}\ge 3/4$. Then
\begin{equation*}
{\mathbb P}\{|z-\theta(\zeta_1,\dots,\zeta_m)|\le \varepsilon\}\ge 1-e^{-m/8}.
\end{equation*}
\end{lemma}

As in \cite{Hei23a, Hei23c} we will use the randomized norm estimation algorithm from \cite{Hei18}. 
 Let $(Q,\mathcal{Q},\varrho)$ be a probability space and let $1\le v< u
\le \infty$. 
For $n\in{\mathbb{N}}$ define  $A_n^{(1)}=(A_{n,\omega}^{(1)})_{\omega\in\Omega}$ by setting 
 for $\omega\in \Omega$ and $f\in L_u(Q,\mathcal{Q},\varrho)$ 
\begin{eqnarray*}
A_{n,\omega}^{(1)}(f)&=&\left(\frac{1}{n}\sum_{i=1}^n|f(\xi_i(\omega_2))|^v\right)^{1/v},
\end{eqnarray*}
where $\xi_i$ are independent $Q$-valued random variables on a probability space $(\Omega,\Sigma,{\mathbb P})$ with distribution $\varrho$. 
The following is essentially  Proposition 6.3 of \cite{Hei18}, for a self-contained proof we refer to \cite{Hei23a}. 
\begin{proposition}
\label{Rpro:4}
Let $1\le v<u\le \infty$.  
Then there is 
a constant $c>0$ such that for  all probability spaces $(Q,\mathcal{Q},\varrho)$, $f\in L_u(Q,\mathcal{Q},\varrho)$, and  $n\in{\mathbb{N}}$ 
\begin{eqnarray*}
{\mathbb E}\,\left|\|f\|_{L_v(Q,\mathcal{Q},\varrho)}-A_{n,\omega}^{(1)}(f)\right| &\le& cn^{\max\left(1/u-1/v,-1/2\right)}\|f\|_{L_u(Q,\mathcal{Q},\varrho)}.
\end{eqnarray*}
\end{proposition}

\section{Algorithms and Upper Bounds for Mean Computation}
\label{sec:3}
Let $1\le p,u\le\infty$. Throughout the paper we use the notation
\begin{equation}
\label{CK5}
\bar{p}=\min(p,2),\quad \bar{u}=\min(u,2).
\end{equation}
We refer to the definition of mean computation $I^{N_1,N_2}$ given in \eqref{C6}.
Expressed in the terminology of \eqref{M7}, we shall study the problem 
\begin{equation*}
\mathcal{P}^{N_1,N_2}=\left( B_{L_p^{N_1}\big(L_u^{N_2}\big)},\K,I^{N_1,N_2},\K,\Lambda^{N_1,N_2}\right),
\end{equation*}
where $\Lambda^{N_1,N_2}$ is standard information consisting of function values, that is,
\begin{equation}
\label{CI7}
\Lambda^{N_1,N_2}=\{\delta_{ij}:\, 1\le i\le N_1,\,1\le j\le N_2\} \text{ with } \delta_{ij}(f)=f(i,j).
\end{equation}
Clearly, this problem is linear. Moreover, we have
\begin{equation}
\label{AF2}
\big\|I^{N_1,N_2}:L_p^{N_1}\big(L_u^{N_2}\big)\to\K\big\|=1.
\end{equation}
We also use the notation $I^{N_2}$ for the mean  operator in $L_u^{N_2}$.
Furthermore, we need the operator of vector valued mean computation $S^{N_1,N_2}$
\begin{equation}
\label{AH1}
S^{N_1,N_2}: L_p^{N_1}\big(L_u^{N_2}\big)\to L_p^{N_1},\quad  (S^{N_1,N_2}f)(i)=\frac{1}{N_2} \sum_{j=1}^{N_2} f(i,j)\quad(i=1,\dots,N_1).
\end{equation}

Given  $n\in\N$, we define for $I^{N_1,N_2}$ a non-adaptive randomized algorithm 
$A_n^{(2)}=\big(A_{n,\omega}^{(2)}\big)_{\omega\in\Omega}$ 
with
$(\Omega,\Sigma,\mu)$ a suitable probability space as follows. Let $\zeta_l \;(l=1,\dots ,n)$ be independent uniformly distributed on 
$\Z[1,N_1]\times\Z[1,N_2]$ random variables, defined on $(\Omega,\Sigma,\mu)$. 
We put for $f\in L_p^{N_1}\big(L_u^{N_2}\big)$,
\begin{eqnarray*}
A_{n,\omega}^{(2)} (f) &=& \frac{1}{n} \sum_{l=1}^{n} f(\zeta_l(\omega)).
\end{eqnarray*}
Note that
the constants in the subsequent statements and proofs are independent of the parameters $n$, $N_1$,$N_2$, and $m$. This is also made clear by the order of quantifiers in the respective statements. 
\begin{proposition}
\label{pro:4}
Let $1 \le p, u\le \infty $ and recall \eqref{CK5}. Then there is a constant $c>0$ such that for all $n,N_1,N_2\in \N$  and $f\in L_p^{N_1}\big(L_u^{N_2}\big)$ 
\begin{equation}
\label{K1}
\E A_{n,\omega}^{(2)}(f)= I^{N_1,N_2}f,\quad \ca(A_{n,\omega}^{(2)})= n,
\end{equation}
and for $p\ge u$ 
\begin{eqnarray}
	\label{K4}
		 \left(\E |I^{N_1,N_2}f-A_{n,\omega}^{(2)}(f)|^{\bar{u}}\right)^{1/\bar{u}}
		&\le &c\min(N_2^{1/\bar{u}-1/\bar{p}}n^{-1+1/\bar{p}},n^{-1+1/\bar{u}}),
\end{eqnarray}
while for $p\le u$
\begin{eqnarray}
	\label{K2}
		 \left(\E |I^{N_1,N_2}f-A_{n,\omega}^{(2)}(f)|^{\bar{p}}\right)^{1/\bar{p}}
		&\le &c\min(N_1^{1/\bar{p}-1/\bar{u}}n^{-1+1/\bar{u}},n^{-1+1/\bar{p}}).
\end{eqnarray}
\end{proposition}
\begin{proof}
Relation \eqref{K1} is obvious. To show \eqref{K4} and \eqref{K2}, we use the factorizations
\begin{eqnarray}
I^{N_1,N_2}&:& L_p^{N_1}\big(L_u^{N_2}\big)\overset{J_1}{\longrightarrow}L_{\bar{p}}^{N_1}\big(L_{\bar{u}}^{N_2}\big)
\overset{J_2}{\longrightarrow}L_{\bar{p}}^{N_1}(L_{\bar{p}}^{N_2})\overset{I^{N_1,N_2}}{\longrightarrow}\K,
\label{CE3}\\
I^{N_1,N_2}&:& L_p^{N_1}\big(L_u^{N_2}\big)\overset{J_1}{\longrightarrow}L_{\bar{p}}^{N_1}\big(L_{\bar{u}}^{N_2}\big)
\overset{J_3}{\longrightarrow}L_{\bar{u}}^{N_1}(L_{\bar{u}}^{N_2})\overset{I^{N_1,N_2}}{\longrightarrow}\K,
\label{CE4}
\end{eqnarray}
with $J_1,J_3$ the identical embeddings.  Then $\|J_1\|=1$, 
\begin{eqnarray}
&&\|J_2\|=N_2^{1/\bar{u}-1/\bar{p}}, \quad\|J_3\|=1\quad \hspace{1.3cm}(p\ge u)\label{CE1}
\\
&&\|J_2\|=1, \quad \hspace{1.3cm}\|J_3\|=N_1^{1/\bar{p}-1/\bar{u}}\quad(p\le u).\label{CE2}
\end{eqnarray}
Furthermore,  the well-known estimate of the Monte Carlo method for 
$L_{w}^{N_1}(L_{w}^{N_2})=L_{w}^{N_1N_2}$ gives
\begin{eqnarray}
\label{CE5}
		\sup_{f\in B_{L_{w}^{N_1N_2}}} \left(\E |I^{N_1,N_2}f-A_{n,\omega}^{(2)}(f)|^w\right)^{1/w}
		&\le &cn^{-1+1/w}\quad(w\in\{\bar{u},\bar{p}\}).
	\end{eqnarray}
Now \eqref{CE3}--\eqref{CE1} and \eqref{CE5} imply \eqref{K4}, while  \eqref{CE3}, \eqref{CE4}, \eqref{CE2}, and \eqref{CE5} yield \eqref{K2}.
 
\end{proof}
\begin{remark} Obviously, for $\bar{p}\ne\bar{u}$ the minimum in \eqref{K4} is attained at the first term iff $n\ge N_2$, and in \eqref{K2} iff $n\ge N_1$.
\end{remark}

Now we define an adaptive randomized algorithm for  the case $1\le p<2<u \le \infty$. Let $f\in L_p^{N_1}\big(L_u^{N_2}\big)$ and set $f_i=(f(i,j))_{j=1}^{N_2}$.  Let $m,n\in\N$, $n\ge N_1$, and let
\begin{equation*}
\left\{\xi_{jk}:\, 1\le j\le \left\lceil\frac{n}{N_1}\right\rceil,\,1\le k\le m\right\}, \quad \left\{\eta_{ij}:\, 1\le i\le N_1,\,1\le j\le n\right\}
\end{equation*}
be independent random variables on a probability space $(\Omega,\Sigma,\Prm)$ uniformly distributed over $\{1,\dots,N_2\}$. We will assume that $(\Omega,\Sigma,\Prm)=(\Omega_1,\Sigma_1,\Prm_1)\times (\Omega_2,\Sigma_2,\Prm_2)$, that the $(\xi_{jk})$ are defined on $(\Omega_1,\Sigma_1,\Prm_1)$, and the $(\eta_{ij})$ on $(\Omega_2,\Sigma_2,\Prm_2)$. Let $\EE,\EE_1,\EE_2$ denote the expectations with respect to the corresponding  probability spaces.

We start by applying  $m$ times algorithm $A_{\left\lceil n/N_1\right\rceil}^{(1)}$ to estimate $\|f_i\|_{L_2^{N_2}}$. Then we compute the median of the results, i.e., we put for $\omega_1\in\Omega_1$, $1\le i\le N_1$,  $1\le k\le m$
\begin{eqnarray*}
a_{ik}(\omega_1)=\Bigg(\left\lceil\frac{n}{N_1}\right\rceil^{-1}\sum_{1\le j\le \left\lceil\frac{n}{N_1}\right\rceil} |f_i(\xi_{jk}(\omega_1))|^2\Bigg)^{1/2},
\quad
\tilde{a}_i(\omega_1)=\theta\big((a_{ik}(\omega_1))_{k=1}^{m}\big).
\end{eqnarray*}
Next we approximate $I^{N_2}f_i$ for each $i$ and $\omega_2\in\Omega_2$ by 
\begin{eqnarray}
b_i(\omega_1,\omega_2)&=&\frac{1}{n_i(\omega_1)}\sum_{j=1}^{n_i(\omega_1)}f_i(\eta_{ij}(\omega_2)). 
\label{UN2}
\end{eqnarray}
where
\begin{equation}
\label{UM5}
 n_i(\omega_1)= \left\{\begin{array}{lll}
\displaystyle \left\lceil\frac{n}{N_1}\right\rceil& \quad\mbox{if}\quad  \tilde{a}_i(\omega_1)^p\le N_1^{-1}\displaystyle \sum_{l=1}^{N_1}\tilde{a}_l(\omega_1)^p  \\[.5cm]
  \displaystyle\left\lceil\frac{\tilde{a}_i^p n}{\sum_{l=1}^{N_1}\tilde{a}_l^p}\right\rceil & \quad\mbox{if}\quad  \tilde{a}_i(\omega_1)^p> N_1^{-1}\displaystyle \sum_{l=1}^{N_1}\tilde{a}_l(\omega_1)^p.    
    \end{array}
\right. 
\end{equation}
Finally we define the output $A_{n,m,\omega}^{(3)}(f)\in \K $ as 
\begin{equation}
\label{UN3}
A_{n,m,\omega}^{(3)}(f)=b(\omega_1,\omega_2):=\frac{1}{N_1}\sum_{i=1}^{N_1}b_i(\omega_1,\omega_2).
\end{equation}
\begin{proposition}
\label{pro:5}
Let $1\le p< 2< u \le  \infty$. Then there exist constants $c_1,c_2>0$ such that the following hold for all $m,n,N_1,N_2\in\N$ with $n\ge N_1$ and $ f\in L_p^{N_1}\big(L_u^{N_2}\big)$:
\begin{equation}
\label{WM1}
\ca(A_{n,m,\omega}^{(3)})\le 6mn
\end{equation}
and for $m\ge c_1\log(N_1+1)$
\begin{eqnarray}
\label{mc-eq:2}
\left(\E |I^{N_1,N_2}f-A_{n,m,\omega}^{(3)}(f)|^{2}\right)^{1/2}
&\le& 
 c_2\Big(N_1^{1/p-1/u}n^{-1+1/u}+n^{-1/2}\Big)\|f\|_{L_p^{N_1}\big(L_u^{N_2}\big)}.\quad
\end{eqnarray}
\end{proposition}
\begin{proof}
The proof is a modification of the  proof of Proposition 4.3 in \cite{Hei23a}. A similar reasoning has been used for Proposition 2 of \cite{Hei23c}.   We have $1\le n_i\le n$ and 
\begin{eqnarray*}
\ca(A_{n,m,\omega}^{(3)})=mN_1\left\lceil\frac{n}{N_1}\right\rceil+\sum_{i=1}^{N_1} n_i
&\le&
 (m+1)N_1\left\lceil\frac{n}{N_1}\right\rceil+
\sum_{i=1}^{N_1}\left\lceil\frac{\tilde{a}_i^p n}{\sum_{l=1}^{N_1}\tilde{a}_l^p}\right\rceil
\nonumber\\[.2cm]
&\le&
 (m+2)(n+N_1)\le 6mn,
\end{eqnarray*}
hence \eqref{WM1}. By definition \eqref{UM5} 
\begin{equation}
\label{UM0}
n_i(\omega_1)\ge \left\lceil\frac{n}{N_1}\right\rceil\quad (1\le i\le N_1).
\end{equation}
By \eqref{UN2}
\begin{eqnarray}
\left(\EE_2\left|\big(S^{N_1,N_2}f\big)(i)-b_i(\omega_1,\omega_2)\right|^2\right)^{1/2}&=&\Bigg(\EE_2\bigg|I^{N_2}f_i-\frac{1}{n_i(\omega_1)}\sum_{j=1}^{n_i(\omega_1)}f_i(\eta_{ij}(\omega_2))\bigg|^2\bigg)^{1/2}
\notag\\
&\le& 
n_i(\omega_1)^{-1/2}\|f_i\|_{L_2^{N_2}} \quad (1\le i\le N_1, \omega_1\in\Omega_1).
\label{UM1}
\end{eqnarray}
Moreover, for fixed $\omega_1\in\Omega_1$  the random variables $b_i(\omega_1,\omega_2)$ are independent and 
\begin{equation*}
\EE_2\, b_i(\omega_1,\omega_2)=\big(S^{N_1,N_2}f\big)(i),
\end{equation*}
therefore 
 by \eqref{UN3}
\begin{eqnarray}
\lefteqn{\left(\EE_2\left|\big(I^{N_1,N_2}f\big)-b(\omega_1,\omega_2)\right|^2\right)^{1/2}}
\notag\\
&=&\left(\EE_2\left|\frac{1}{N_1}\sum_{i=1}^{N_1}\big(\big(S^{N_1,N_2}f\big)(i)-b_i(\omega_1,\omega_2)\Big)\right|^2\right)^{1/2}
\nonumber\\
&=& N_1^{-1/2}\left(\frac{1}{N_1}\sum_{i=1}^{N_1}\EE_2\left|\big(\big(S^{N_1,N_2}f\big)(i)-b_i(\omega_1,\omega_2)\Big)\right|^2\right)^{1/2}.\quad (\omega_1\in\Omega_1)
\label{AD4}
\end{eqnarray}
We set
$
c_1=16/\log e,
$
and assume $m\ge c_1\log(N_1+1)$, hence $e^{-m/8}\le (N_1+1)^{-2}$.
Let $c(2)$ denote the constant from Proposition \ref{Rpro:4}, which implies  for fixed $f\in L_p^{N_1}\big(L_u^{N_2}\big)$
\begin{equation*}
\Prm_1\left\{\omega_1\in \Omega_1:\,\Big|\|f_i\|_{L_2^{N_2}}-a_{ik}(\omega_1)\Big|\le 4c(2) \bigg(\frac{n}{N_1}\bigg)^{-\left(1/2-1/u\right)}\|f_i\|_{L_u^{N_2}}\right\}\ge \frac{3}{4},
\end{equation*}
so by Lemma \ref{Ulem:2e}
\begin{equation*}
\Prm_1\left\{\omega_1\in \Omega_1:\,\Big|\|f_i\|_{L_2^{N_2}}-\tilde{a}_i(\omega_1)\Big|\le 4c(2) \bigg(\frac{n}{N_1}\bigg)^{-\left(1/2-1/u\right)}\|f_i\|_{L_u^{N_2}}\right\}\ge 1-(N_1+1)^{-2}.
\end{equation*}
Setting 
\begin{equation}
\label{L0}
\Omega_{1,0}=\left\{\omega_1\in \Omega_1:\,\Big|\|f_i\|_{L_2^{N_2}}-\tilde{a}_i(\omega_1)\Big|\le 4c(2) \bigg(\frac{n}{N_1}\bigg)^{-\left(1/2-1/u\right)}\|f_i\|_{L_u^{N_2}}\quad(1\le i\le N_1)\right\},
\end{equation}
we arrive at
\begin{equation}
\label{L1}
\Prm_1(\Omega_{1,0})\ge 1-(N_1+1)^{-1}.
\end{equation}
 Fix $\omega_1\in \Omega_{1,0}$.
Then by \eqref{L0} for all $i$
\begin{equation*}
\tilde{a}_i(\omega_1)\le c \|f_i\|_{L_u^{N_2}}.
\end{equation*}
Consequently,
\begin{equation}
\label{UL9}
\frac{1}{N_1}\sum_{i=1}^{N_1}\tilde{a}_i(\omega_1)^p\le \frac{c}{N_1}\sum_{i=1}^{N_1}\|f_i\|_{L_u^{N_2}}^p
= c\|f\|_{L_p^{N_1}\big(L_u^{N_2}\big)}^p.
\end{equation}
Set 
\begin{equation*}
I(\omega_1):=\bigg\{1\le i\le N_1:\, \tilde{a}_i(\omega_1)\le \frac{\|f_i\|_{L_2^{N_2}}}{2}\bigg\},
\end{equation*}
then from \eqref{L0}
\begin{equation*}
\|f_i\|_{L_2^{N_2}}\le c\bigg(\frac{n}{N_1}\bigg)^{-\left(1/2-1/u\right)}\|f_i\|_{L_u^{N_2}} \quad (i\in I(\omega_1)),
\end{equation*}
which together with \eqref{UM0} and \eqref{UM1} gives 
\begin{eqnarray}
\left(\EE_2\left|\big(S^{N_1,N_2}f\big)(i)-b_i(\omega_1,\omega_2)\right|^2\right)^{1/2}\le 
c\bigg(\frac{n}{N_1}\bigg)^{-1+1/u}\|f_i\|_{L_u^{N_2}}\quad (i\in I(\omega_1)) .
\label{UM2}
\end{eqnarray}
Now let $i\not\in I(\omega_1)$, thus 
\begin{equation}
\label{UM6}
\tilde{a}_i(\omega_1)> \frac{\|f_i\|_{L_2^{N_2}}}{2}.
\end{equation}
If $\tilde{a}_i^p\le N_1^{-1}\sum_{l=1}^{N_1}\tilde{a}_l^p$, then by \eqref{UM5},  \eqref{UM1}, \eqref{UL9}, and \eqref{UM6}
\begin{eqnarray}
\left(\EE_2\left|\big(S^{N_1,N_2}f\big)(i)-b_i(\omega_1,\omega_2)\right|^2\right)^{1/2}&\le& 
2\bigg(\frac{n}{N_1}\bigg)^{-1/2}\tilde{a}_i
\le 2\bigg(\frac{n}{N_1}\bigg)^{-1/2}\bigg(\frac{\sum_{l=1}^{N_1}\tilde{a}_l^p}{N_1}\bigg)^{1/p}
\nonumber\\
&\le& c\bigg(\frac{n}{N_1}\bigg)^{-1/2}\|f\|_{L_p^{N_1}\big(L_u^{N_2}\big)}.
\label{UM4}
\end{eqnarray}
Similarly, if $\tilde{a}_i^p> N_1^{-1} \sum_{l=1}^{N_1}\tilde{a}_l^p$, then 
\begin{eqnarray}
\lefteqn{\left(\EE_2\left|\big(S^{N_1,N_2}f\big)(i)-b_i(\omega_1,\omega_2)\right|^2\right)^{1/2}\le 
2\bigg(\frac{\tilde{a}_i^p n}{\sum_{l=1}^{N_1}\tilde{a}_l^p}\bigg)^{-1/2}\tilde{a}_i}
\nonumber\\
&=&2\tilde{a}_i^{1-p/2}\bigg(\frac{n}{N_1}\bigg)^{-1/2}\bigg(\frac{\sum_{l=1}^{N_1}\tilde{a}_l^p}{N_1}\bigg)^{1/2}
\le c\tilde{a}_i^{1-p/2}\bigg(\frac{n}{N_1}\bigg)^{-1/2}\|f\|_{L_p^{N_1}\big(L_u^{N_2}\big)}^{p/2}.
\label{U5}
\end{eqnarray}
From \eqref{AD4}, \eqref{UL9}, \eqref{UM2}, \eqref{UM4}, and \eqref{U5} we obtain
\begin{eqnarray*}
\lefteqn{\left(\EE_2\left|\big(I^{N_1,N_2}f\big)-b(\omega_1,\omega_2)\right|^2\right)^{1/2}}
\notag\\
&=& N_1^{-1/2}\left(\frac{1}{N_1}\sum_{i=1}^{N_1}\EE_2\left|\big(\big(S^{N_1,N_2}f\big)(i)-b_i(\omega_1,\omega_2)\Big)\right|^2\right)^{1/2}
\nonumber\\
&\le&cN_1^{-1/2}\bigg(\frac{n}{N_1}\bigg)^{-1+1/u}\left(\frac{1}{N_1}\sum_{i=1}^{N_1}\|f_i\|_{L_u^{N_2}}^2\right)^{1/2}+cn^{-1/2}\|f\|_{L_p^{N_1}\big(L_u^{N_2}\big)}
\nonumber\\
&&+c\left(\frac{1}{N_1}\sum_{i=1}^{N_1}\tilde{a}_i^{2-p}\right)^{1/2}n^{-1/2}\|f\|_{L_p^{N_1}\big(L_u^{N_2}\big)}^{p/2}.
\label{UN4}
\end{eqnarray*}
Observe that by H\"older's inequality
\begin{eqnarray*}
\left(\frac{1}{N_1}\sum_{i=1}^{N_1}\|f_i\|_{L_u^{N_2}}^2\right)^{1/2}&\le&N_1^{1/p-1/2}\left(\frac{1}{N_1}\sum_{i=1}^{N_1}\|f_i\|_{L_u^{N_2}}^p\right)^{1/p} = N_1^{1/p-1/2} \|f\|_{L_p^{N_1}\big(L_u^{N_2}\big)} 
\\
\left(\frac{1}{N_1}\sum_{i=1}^{N_1}\tilde{a}_i^{2-p}\right)^{1/2}&\le& \left(\frac{1}{N_1}\sum_{i=1}^{N_1}\tilde{a}_i^p\right)^{\frac{2-p}{2p}}\le \|f\|_{L_p^{N_1}\big(L_u^{N_2}\big)}^{1-p/2}.
\end{eqnarray*}
Inserting this into \eqref{UN4} gives 
\begin{eqnarray}
\lefteqn{\left(\EE_2\left|\big(I^{N_1,N_2}f\big)-b(\omega_1,\omega_2)\right|^2\right)^{1/2}}
\nonumber\\
&\le&c\Big(N_1^{1/p-1/u}n^{-1+1/u}+n^{-1/2}\Big)\|f\|_{L_p^{N_1}\big(L_u^{N_2}\big)} \quad(\omega_1\in\Omega_{1,0}).
\label{UN1}
\end{eqnarray}

Finally we estimate the error on $(\Omega_1\setminus \Omega_{1,0})\times\Omega_2$. For $\omega_1\in \Omega_1$ we conclude from \eqref{UM0}--\eqref{AD4}
\begin{eqnarray*}
\lefteqn{\left(\EE_2\left|\big(I^{N_1,N_2}f\big)-b(\omega_1,\omega_2)\right|^2\right)^{1/2}}
\notag\\
&\le& N_1^{-1/2}\bigg(\frac{n}{N_1}\bigg)^{-1/2}\left(\frac{1}{N_1}\sum_{i=1}^{N_1}\|f_i\|_{L_2^{N_2}}^2\right)^{1/2}\le n^{-1/2}\left(\frac{1}{N_1}\sum_{i=1}^{N_1}\|f_i\|_{L_u^{N_2}}^2\right)^{1/2}
\notag\\
&\le&
N_1^{1/p-1/2} n^{-1/2}\left(\frac{1}{N_1}\sum_{i=1}^{N_1}\|f_i\|_{L_u^{N_2}}^p\right)^{1/p} = N_1^{1/p-1/2}  n^{-1/2}\|f\|_{L_p^{N_1}\big(L_u^{N_2}\big)}. 
\end{eqnarray*}
Consequently, using also \eqref{L1},
\begin{eqnarray*}
\lefteqn{\left(\int_{(\Omega_1\setminus \Omega_{1,0})\times\Omega_2}\big|I^{N_1,N_2}f-b(\omega)\big|^2 d\Prm(\omega)\right)^{1/2}}
\\
&\le& \Prm_1(\Omega_1\setminus \Omega_{1,0})^{1/2}N_1^{1/p-1/2}  n^{-1/2}\|f\|_{L_p^{N_1}\big(L_u^{N_2}\big)}
\\
&\le& (N_1+1)^{-1/2} N_1^{1/p-1/2} n^{-1/2}\|f\|_{L_p^{N_1}\big(L_u^{N_2}\big)}<n^{-1/2}\|f\|_{L_p^{N_1}\big(L_u^{N_2}\big)},
\end{eqnarray*}
which together with \eqref{UN1} proves \eqref{mc-eq:2}.

\end{proof}
\section{Lower Bounds and Complexity} \label{sec:4}

\begin{proposition}
\label{pro:1}
Let $1\le p,u\le \infty$. Then there are constants $0<c_0<1$, $c_{1-4}>0$ such that for each $N_1,N_2\in \N$ there exist probability measures $\mu_{N_1,N_2}^{(1)},\dots,\mu_{N_1,N_2}^{(4)}$ with finite support in $B_{L_p^{N_1}\big(L_u^{N_2}\big)}$ such that for  $n\in \N$ with $n<c_0N_1N_2$
\begin{eqnarray}
e_n^\avg(I^{N_1,N_2},\mu_{N_1,N_2}^{(1)})
&\ge&  c_1N_1^{1/p-1}N_2^{1/u-1}
\label{UL2}\\
e_n^\avg(I^{N_1,N_2},\mu_{N_1,N_2}^{(2)})
&\ge& c_2 (N_1N_2)^{-1/2}
\label{UL1}\\
e_n^\avg(I^{N_1,N_2},\mu_{N_1,N_2}^{(3)})
&\ge&  c_3N_1^{-1/2}N_2^{1/u-1}
\label{CB2}\\
e_n^\avgno(I^{N_1,N_2},\mu_{N_1,N_2}^{(4)})
&\ge& c_4 N_1^{1/p-1}N_2^{-1/2}.\label{UQ3}
\end{eqnarray}
\end{proposition}
\begin{proof}
We set $c_0=\frac1{21}$ and let $n\in \N$ be such that 
\begin{equation}
\label{B11}
n< \frac{N_1N_2}{21}.
\end{equation}
To prove \eqref{UL2}, we define functions $\psi_{ij}\in L_p^{N_1}\big(L_u^{N_2}\big)$ by
\begin{equation*}
\psi_{ij}(s,t) = N_1^{1/p}N_2^{1/u}e_{ij} \quad(i\in\Z[1,N_1],\,j\in \Z[1,N_2]),
\end{equation*}
with $(e_{ij})$ the unit vectors in $\K^{N_1\times N_2}$,
and let $\mu_{N_1,N_2}^{(1)}$ be the uniform distribution on the set 
$$
\{ \alpha \psi_{ij}:\,i=1, \dots, N_1,\, j=1, \dots, N_2 ,\, \alpha =\pm 1\} \subset B_{L_p^{N_1}\big(L_u^{N_2}\big)} .
$$
Recall that by (\ref{B11}), $4n<N_1N_2$, so from Lemma \ref{lem:5}(ii) a we conclude
\begin{eqnarray*}
e_n^\avg(I^{N_1,N_2},\mu_{N_1,N_2}^{(1)})
&\ge& \frac12\big| I^{N_1,N_2} \psi_{1,1}\big|
 =  \frac12 N_1^{1/p-1}N_2^{-1+1/u},
\end{eqnarray*}
thus \eqref{UL2}.

To show the second lower bound, \eqref{UL1}, 
  let $(\varepsilon_{ij})_{i=1,j=1}^{N_1,N_2}$ be independent 
symmetric Bernoulli random variables and let $\mu_{N_1,N_2}^{(2)}$ be the distribution of 
$
\sum_{i=1,j=1}^{N_1,N_2}\epsilon_{ij}e_{ij}.
$
Since by (\ref{B11}), $4n<N_1N_2$, we can apply Lemma \ref{lem:5}. So let $\mathcal{K}$ be any subset of 
$\{(i,j)\,:\,1\le i\le N_1,\, 1\le j\le N_2\}$ with $|\mathcal{K}|\ge N_1N_2-2n$. Then $|\mathcal{K}|\ge N_1N_2/2$
and we obtain
from   Khintchine's inequality, 
\begin{eqnarray*}
\E \bigg|   \sum_{(i,j)\in \mathcal{K}} \epsilon_{ij} I^{N_1,N_2}e_{ij}  \bigg|
 =\frac{1}{N_1N_2}\E \bigg|  \sum_{(i,j)\in \mathcal{K}} \varepsilon_{ij} \bigg|
  \ge \frac{c|\mathcal{K}|^{1/2}}{N_1N_2}\ge c(N_1N_2)^{-1/2},
\end{eqnarray*}
and therefore from Lemma \ref{lem:5} (i),
\begin{eqnarray*}
e_n^\avg(I^{N_1,N_2},\mu_{N_1,N_2}^{(2)})&\ge& \frac{1}{2}\, \min_{|\mathcal{K}|\ge N_1N_2-2n} \E\bigg|   \sum_{(i,j)\in \mathcal{K}} \varepsilon_{ij} I^{N_1,N_2}\psi_{ij}  \bigg|
 \ge c(N_1N_2)^{-1/2}.
\end{eqnarray*}
This proves \eqref{UL1}.

Next we turn to \eqref{UQ3}, where we use Corollary \ref{Ucor:2} with 
\begin{eqnarray}
 &&M=N_1, \quad F_1=L_u^{N_2},\quad S_1=I^{N_2},\quad G_1=K_1=\K,\quad \Lambda_1=\{\delta_j:\, 1\le j\le N_2\}, 
\label{L5}\\
&& f_{i,0}'=0\quad(i=1,\dots,N_1),
\label{AB0}
\end{eqnarray}
where $\delta_j(g)=g(j)$. Then obviously \eqref{J1} is satisfied and
\begin{eqnarray}
&&F=\prod_{i=1}^{N_1} L_u=  L_p^{N_1}\big(L_u^{N_2}\big),\quad G=\K,\quad S=N_1 I^{N_1,N_2},
\label{L4}\\
&& K=\K,\quad \Lambda=\Lambda^{N_1,N_2}=\{\delta_{ij}:\, 1\le i\le N_1,\, 1\le j\le N_2\}. 
\label{L8}
\end{eqnarray}
We define 
\begin{equation}
\label{UP3}
\psi_j=N_1^{1/p}e_j\in L_u^{N_2} \quad (j=1,\dots,N_2).
\end{equation}
Let $(\epsilon_j)_{j=1}^{N_2}$ be independent symmetric Bernoulli random variables, let $\mu_1$ be the distribution of  $ \sum_{j=1}^{N_2}  \epsilon_j\psi_j$, and denote the resulting from \eqref{B8} measure by $\mu_{N_1,N_2}^{(4)}$. Observe that by \eqref{UP3} $\mu_{N_1,N_2}^{(4)}$ is supported by $B_{L_p^{N_1}\big(L_u^{N_2}\big)}$. 
Now  \eqref{U0} and \eqref{B7} yield
\begin{eqnarray}
\label{UP2}
e_n^\avgno(I^{N_1,N_2},\mu_{N_1,N_2}^{(4)})&=& N^{-1}e_n^\avgno(S,\mu_{N_1,N_2}^{(4)})\ge\frac{1}{2} N_1^{-1}e_{\left\lfloor\frac{2n}{N_1}\right\rfloor}^\avgno(I^{N_2},\mu_1,\K)
\nonumber\\
&\ge& \frac{1}{2} N_1^{-1}e_{\left\lfloor\frac{2n}{N_1}\right\rfloor}^\avg(I^{N_2},\mu_1,\K).
\end{eqnarray}
By \eqref{B11} we have $4\left\lfloor\frac{2n}{N_1}\right\rfloor<N_2$, therefore Lemma \ref{lem:5}(i) and Khintchine's inequality give
\begin{eqnarray*}
e_{\left\lfloor\frac{2n}{N_1}\right\rfloor}^\avg(I^{N_2},\mu_1,\K)&\ge& \frac{1}{2}\min\bigg\{\E\Big|\sum_{i\in \mathcal{I}}\epsilon_iI^{N_2}\psi_i\Big|:\,\mathcal{I}\subseteq\{1,\dots,N_2\},\,|\mathcal{I}|\ge N_2-2\left\lfloor\frac{2n}{N_1}\right\rfloor\bigg\}
\\
&\ge& c N_2^{1/2}|I^{N_2}\psi_1|\ge cN_1^{1/p}N_2^{-1/2}.
\end{eqnarray*}
Inserting this into \eqref{UP2}  yields \eqref{UQ3}.

Finally we prove \eqref{CB2}. For this purpose we can assume $N_2>16$, because for $N_2\le 16$ relation \eqref{CB2} follows from \eqref{UL1}. 
We define 
\begin{equation}
\label{CB3}
\psi_{j}(t) = \left\{ \begin{array}{ll}
                N_2^{1/u} & \text{ if} \quad t=j,\\
                0  & \text{ otherwise.}
                         \end{array}
                 \right.
\end{equation}
Let $\mu_{1}$ be the uniform distribution on $\{\pm \psi_j:\, 1\le j\le N_2\}\subset L_u^{N_2}$. Then we set $\mu_{N_1,N_2}^{(4)}=\mu_1^{N_1}$. 
Due to \eqref{CB3} this measure has its support in $B_{L_p^M(L_u^{N_2})}$. First we need an auxiliary estimate for the operator of vector-valued mean computation $S^{N_1,N_2}:L_p^{N_1}(L_u^{N_2})\to L_1^{N_1}$, see \eqref{AH1}. With the choice \eqref{L5}--\eqref{AB0} we are in the setting of 
 Corollary 2.4 of \cite{Hei23a} and \eqref{L4} and \eqref{L8} hold, except that now $G=L_1^{N_1}$ and $S=S^{N_1,N_2}$.  It follows that
\begin{equation}
\label{B9}
e_n^{\rm avg }(S^{N_1,N_2},\mu_{N_1,N_2}^{(4)},L_1^{N_1})\ge \frac14 e_{\left\lceil\frac{4n}{N_1}\right\rceil}^{\rm avg }(I^{N_2},\mu_1).
\end{equation}
We have 
\begin{equation}
\label{CB4}
4\left\lceil\frac{4n}{N_1}\right\rceil+1\le N_2.
\end{equation}
Indeed, this is clear for $n\le N_1$, since we assumed $N_2>16$, while for $n> N_1$ assumption \eqref{B11} implies   
\begin{equation*}
4\left\lceil\frac{4n}{N_1}\right\rceil+1\le \frac{16n}{N_1}+5\le \frac{21n}{N_1}\le N_2.
\end{equation*}
By Lemma \ref{lem:5}(ii) with $\bar{n}=N_2$, taking into account \eqref{CB4},
\begin{eqnarray*}
e_{\left\lceil\frac{4n}{N_1}\right\rceil}^{\rm avg }(I^{N_2},\mu_1)\ge \frac{1}{2}|I^{N_2}\psi_1|=\frac{1}{2}N_2^{1/u-1},
\end{eqnarray*}
which together with  \eqref{B9} gives 
\begin{equation}
\label{CB5}
e_n^{\rm avg }(S^{N_1,N_2},\mu_{N_1,N_2}^{(4)},L_1^{N_1})\ge \frac18 N_2^{1/u-1}.
\end{equation}

Now let $A=((L_k)_{k=1}^\infty, (\tau_k)_{k=0}^\infty,(\varphi_k)_{k=0}^\infty)$  be a deterministic algorithm for (scalar-valued) mean computation $I^{N_1,N_2}$ with 
\begin{equation}
\label{AB2}
\ca(A,\mu_{N_1,N_2}^{(4)})\le n.
\end{equation}
Let $f\in \SUPP( \mu_{N_1,N_2}^{(4)})$, which means that 
\begin{equation*}
f=\sum_{i=1}^{N_1}\alpha_f(i) N_2^{1/u}e_{i,j_f(i)},
\end{equation*}
where $\alpha_f:\Z[1,N_1]\to\{-1,+1\}$ and $j_f:\Z[1,N_1]\to \Z[1,N_2]$. Let $\lambda_{1,f}, \dots, \lambda_{n(f),f}$ be the sequence of information functionals called by $A$ at input $f$ and define 
\begin{equation*}
I_f=\big\{i\in \Z[1,N_1]:\, \delta_{i,j_f(i)}\in \{\lambda_{1,f}, \dots, \lambda_{n(f),f}\} \big\}.
\end{equation*}

Based on $A$ we define an algorithm $\tilde{A}=((L_k)_{k=1}^\infty, (\tau_k)_{k=0}^\infty,(\tilde{\phi}_k)_{k=0}^\infty)$ for vector valued mean computation $S^{N_1,N_2}:L_p^{N_1}(L_u^{N_2})\to L_1^{N_1}$ considered above, that is, the information and stopping mappings $(L_k)$ and $(\tau_k)$ are the same, just the output mappings $\tilde{\phi}_k:\K^k\to L_1^{N_1}$ are defined differently as follows. 
Let  $\tilde{\phi}_0=0\in L_1^{N_1}$. For $k\ge 1$ let $(a_1,\dots, a_k)\in \K^k$ be given, let
\begin{equation*}
\lambda_1=L_1, \quad\lambda_2=L_2(a_1),\quad\dots, \quad\lambda_k=L_k(a_1,\dots, a_{k-1}),
\end{equation*}
and assume that $\tilde{\phi}_{k-1}(a_1,\dots,a_{k-1})\in L_1^{N_1}$  has already been defined. Let $i,j$ be such that $\lambda_k=\delta_{ij}$ . Then we set
\begin{equation*}
\tilde{\phi}_k(a_1,\dots,a_k)=\left\{\begin{array}{lll}
   \tilde{\phi}_{k-1}(a_1,\dots,a_{k-1})+N_2^{-1}a_k e_i &\quad\mbox{if}\quad \lambda_k\not\in \{\lambda_1, \dots,\lambda_{k-1}\}   \\
   \tilde{\phi}_{k-1}(a_1,\dots,a_{k-1})& \quad\mbox{if}\quad\lambda_k\in \{\lambda_1, \dots,\lambda_{k-1}\}. 
    \end{array}
\right. 
\end{equation*}
It is readily checked by induction that for $f\in \SUPP( \mu_{N_1,N_2}^{(4)})$
\begin{equation*}
 (\tilde{A}(f))(i) =\left\{\begin{array}{lll}
 N_2^{-1}f(i,j_f(i))& \quad\mbox{if}\quad i\in I_f   \\
 0  & \quad\mbox{if}\quad i\not\in I_f,
    \end{array}
\right. 
\end{equation*}
while for all $i$
\begin{equation*}
(S^{N_1,N_2}f)(i)=N_2^{-1}f(i,j_f(i)),
\end{equation*}
and therefore
\begin{equation}
\label{CB9}
\|S^{N_1,N_2}f-\tilde{A}(f)\|_{L_1^{N_2}}=N_1^{-1}N_2^{1/u-1}(N_1-|I_f|).
\end{equation}
Clearly, we have 
\begin{equation*}
\ca(\tilde{A},f)=\ca(A,f) \quad (f\in L_p^{N_1}(L_u^{N_2})).
\end{equation*}
Therefore, combining \eqref{CB5}, \eqref{AB2}, and \eqref{CB9} we conclude
\begin{equation}
\label{CC1}
\int_{L_p^{N_1}(L_u^{N_2})} (N_1-|I_f|)d\mu_{N_1,N_2}^{(4)}(f)\ge \frac{N_1}{8}.
\end{equation}
Let 
\begin{equation*}
F_0=\bigg\{f\in \SUPP(\mu_{N_1,N_2}^{(4)}):\ N_1-|I_f|\ge \frac{N_1}{16}\bigg\},
\end{equation*}
then \eqref{CC1} implies 
\begin{equation}
\label{CC0}
\mu_{N_1,N_2}^{(4)}(F_0)\ge \frac1{16}.
\end{equation}
We introduce an equivalence relation on $\SUPP( \mu_{N_1,N_2}^{(4)})$ as follows: $f\sim g$ iff they are non-zero at the same places and their values coincide in all rows $i\in I_f$, or equivalently,
\begin{equation}
f\sim g \quad \text{ iff} \quad j_f=j_g \; \text{and}\; f(i,j_f(i))=g(i,j_f(i))\quad(i\in I_f).\label{CC7}
\end{equation}
It is readily checked by induction that for $f\sim g$ the action of algorithm $A$ is the same on both $f$ and $g$, meaning that for $1\le k\le n(f)$
\begin{equation*}
\lambda_{k,f}=\lambda_{k,g},\quad \lambda_{k,f}(f)=\lambda_{k,f}(g),
\end{equation*}
and therefore also 
\begin{equation*}
n(f)=n(g), \quad I_f=I_g, \quad A(f)=A(g),
\end{equation*}
which shows, in particular, that $\sim$ is indeed an equivalence relation. Moreover, $f\sim g$ and  $f\in F_0$ implies $g\in F_0$.
Let $[f]$ denote the equivalence class of $f\in F_0$ and $[F_0]$ the set of all equivalence classes of elements of $F_0$. Now we estimate
\begin{eqnarray}
\lefteqn{\int_{L_p^{N_1}(L_u^{N_2})} \big|I^{N_1,N_2}f-A(f)\big|d\mu_{N_1,N_2}^{(4)}(f)}
\nonumber\\
&\ge&\sum_{f\in F_0}\big|I^{N_1,N_2}f-A(f)\big|\mu_{N_1,N_2}^{(4)}(\{f\})= (2N_2)^{-N_1}\sum_{[f]\in[F_0]}\sum_{g\in [f]}\big|I^{N_1,N_2}g-A(f)\big|
\nonumber\\
&=& (2N_2)^{-N_1}\sum_{[f]\in[F_0]}\sum_{g\in [f]}\bigg|(N_1N_2)^{-1}\sum_{i\in\Z[1,N_1]\setminus I_f}g(i,j_f(i))+
(N_1N_2)^{-1}\sum_{i\in I_f}f(i,j_f(i))-A(f)\bigg|
\nonumber\\
&=& (2N_2)^{-N_1}\sum_{[f]\in[F_0]}\sum_{g\in [f]}\bigg|-(N_1N_2)^{-1}\sum_{i\in\Z[1,N_1]\setminus I_f}g(i,j_f(i))+
(N_1N_2)^{-1}\sum_{i\in I_f}f(i,j_f(i))-A(f)\bigg|
\nonumber\\
&\ge& (N_1N_2)^{-1}(2N_2)^{-N_1}\sum_{[f]\in[F_0]}\sum_{g\in [f]}\bigg|\sum_{i\in\Z[1,N_1]\setminus I_f}g(i,j_f(i))\bigg|
\label{CC5}\\
&=& N_1^{-1}N_2^{1/u-1}(2N_2)^{-N_1}\sum_{[f]\in[F_0]}\sum_{g\in [f]}\bigg|\sum_{i\in\Z[1,N_1]\setminus I_f}\alpha_g(i)\bigg|,
\nonumber
\end{eqnarray}
where we used a standard symmetry argument to reach \eqref{CC5}. Next observe that according to \eqref{CC7} if $g$ runs through $[f]$, then $(\alpha_g(i))_{i\in\Z[1,N_1]\setminus I_f}$ runs through all possible combinations of $+1$ and $-1$. Therefore Khintchine's inequality gives
\begin{eqnarray*}
\lefteqn{\int_{L_p^{N_1}(L_u^{N_2})} \big|I^{N_1,N_2}f-A(f)\big|d\mu_{N_1,N_2}^{(4)}(f)}
\nonumber\\
&\ge& cN_1^{-1}N_2^{1/u-1}(2N_2)^{-N_1}\sum_{[f]\in[F_0]} |[f]|\big|\Z[1,N_1]\setminus I_f\big|^{1/2}
\nonumber\\
&\ge& cN_1^{-1/2}N_2^{1/u-1}(2N_2)^{-N_1}\sum_{[f]\in[F_0]} |[f]|=cN_1^{-1/2}N_2^{1/u-1}(2N_2)^{-N_1}|F_0| 
\nonumber\\
&=& c\mu_{N_1,N_2}^{(4)}(F_0)N_1^{-1/2}N_2^{1/u-1}\ge cN_1^{-1/2}N_2^{1/u-1},
\end{eqnarray*}
where we used \eqref{CC0} in the last relation. This proves \eqref{CB2}.

\end{proof}
\begin{lemma}
\label{lem:1}
Let $1\le M_\iota\le N_\iota$ $(\iota=1,2)$. Then for $\set\in \{\de,\deno,\ran, \ranno\}$
\begin{equation*}
e_n^\set(I^{M_1,M_2},B_{L_p^{M_1}(L_u^{M_2})})\le 4e_n^\set(I^{N_1,N_2},B_{L_p^{N_1}(L_u^{N_2})}).
\end{equation*}
\end{lemma}
\begin{proof}
We define $M_\iota$ disjoint blocks of $\{1,\dots,N_\iota\}$ by setting
\begin{equation*}
D_{\iota,k}= \left\{ (k-1)\left\lfloor \frac{N_\iota}{M_\iota}\right\rfloor+1, \dots, 
k \left\lfloor \frac{N_\iota}{M_\iota}\right\rfloor  \right\}  \quad (k=1, \dots, M_\iota).
\end{equation*}
We have
\begin{equation*}
\frac{N_\iota}{2M_\iota}\le |D_{\iota,k}|=|D_{\iota,1}|=\left\lfloor \frac{N_\iota}{M_\iota}\right\rfloor \le  \frac{N_\iota}{M_\iota}.
\end{equation*}
With
\begin{equation}
\label{CC4}
\gamma_\iota=\frac{N_\iota}{M_\iota|D_{\iota,1}|}
\end{equation}
it follows that 
\begin{equation}
\label{CC8}
1\le \gamma_\iota\le 2.
\end{equation}
Now we let
$\psi_{kl}\in L_p^{N_1}(L_u^{N_2})$ be defined for $1\le k\le M_1$ and $1\le l\le M_2$ as
\begin{equation*}
\psi_{kl}(s,t) = \left\{ \begin{array}{ll}
                1 & \text{ if} \enspace s\in D_{1,k}\quad\mbox{and}\quad t\in D_{2,l},\\
                0  & \text{ otherwise}
                         \end{array}
                 \right.
\end{equation*}
and  $R:L_p^{M_1}(L_u^{M_2})\to  L_p^{N_1}(L_u^{N_2})$ by
\begin{equation}
Rf=\gamma_1^{1/p}\gamma_2^{1/u}\sum_{k=1}^{M_1}\sum_{l=1}^{M_2}f(k,l)\psi_{kl}.
\label{CK6}
\end{equation}
It follows from \eqref{CC4} that for $f\in L_p^{M_1}(L_u^{M_2})$
\begin{eqnarray}
\|Rf\|_{L_p^{N_1}(L_u^{N_2})}&=&\gamma_1^{1/p}\gamma_2^{1/u}
\left(\frac{1}{N_1}\sum_{k=1}^{M_1}\left(\frac{1}{N_2}\sum_{l=1}^{M_2}|f(k,l)|^u|D_{2,l}|\right)^{p/u}|D_{1,k}|\right)^{1/p}
\notag\\
&=&
\left(\frac{1}{M_1}\sum_{k=1}^{M_1}\left(\frac{1}{M_2}\sum_{l=1}^{M_2}|f(k,l)|^u\right)^{p/u}\right)^{1/p}
=\|f\|_{L_p^{M_1}(L_u^{M_2})}\label{CC9}
\end{eqnarray}
(with the obvious modifications if $p=\infty$ and/or $u=\infty$)
and
\begin{eqnarray}
I^{N_1,N_2}Rf&=&\gamma_1^{1/p}\gamma_2^{1/u}\frac{1}{N_1} \sum_{k=1}^{M_1}\frac{|D_{1,k}|}{N_2}\sum_{l=1}^{M_2}f(k,l)|D_{2,l}|
\nonumber\\
&=&\gamma_1^{1/p}\gamma_2^{1/u}N_1^{-1}N_2^{-1} |D_{1,1}| |D_{2,1}|M_1M_2\,I^{M_1,M_2}f
=\gamma_1^{1/p-1} \gamma_2^{1/u-1}I^{M_1,M_2}f.
\label{CC6}
\end{eqnarray}
Next we show that 
$$
\Big(B_{L_p^{M_1}\big(L_u^{M_2}\big)},\K,I^{M_1,M_2},\K,\Lambda^{M_1,M_2}\Big)\;\text{  reduces to } \; \Big(B_{L_p^{N_1}\big(L_u^{N_2}\big)},\K,I^{N_1,N_2},\K,\Lambda^{N_1,N_2}\Big),
$$ 
so that we can apply Proposition \ref{pro:2} (with $\kappa=1$). We define two mappings 
$$
\eta:\Lambda^{N_1,N_2}\to \Lambda^{M_1,M_2}, \quad\rho:\Lambda^{N_1,N_2}\times \K\to \K
$$
as follows. Let $(i,j)\in\Z[1,N_1]\times \Z[1,N_2]$ and $a\in \K$. If there is a $(k,l)\in\Z[1,M_1]\times \Z[1,M_2]$  such that  $(i,j)\in D_{1,k}\times D_{2,l}$, then we put $\eta(\delta_{ij})=\delta_{kl}$ and $\rho(\delta_{ij},a)=\gamma_1^{1/p}\gamma_2^{1/u}a$. Hence \eqref{CK6} gives for $f\in  L_p^{M_1}(L_u^{M_2})$
\begin{equation}
\label{CK7}
\delta_{ij}(Rf)=\gamma_1^{1/p}\gamma_2^{1/u}\delta_{kl}(f)=\rho(\delta_{ij},\delta_{kl}(f))
=\rho(\delta_{ij},(\eta(\delta_{ij}))(f)).
\end{equation}
If there is no $(k,l)\in\Z[1,M_1]\times \Z[1,M_2]$  with  $(i,j)\in D_{1,k}\times D_{2,l}$ then we
set $\eta(\delta_{ij})=\delta_{1,1}$ and $\rho(\delta_{ij},a)=0$. Here \eqref{CK6} implies that for $f\in L_p^{M_1}(L_u^{M_2})$ 
\begin{equation*}
\delta_{ij}(R(f))=0=\rho(\delta_{ij}, \delta_{1,1}(f))=\rho(\delta_{ij},\delta_{kl}(f)),
\end{equation*}
which together with \eqref{CK7} shows that \eqref{XG3} is satisfied. Furthermore, by \eqref{CC9}
$$
R\Big(B_{L_p^{M_1}\big(L_u^{M_2}\big)}\Big)\subseteq B_{L_p^{N_1}\big(L_u^{N_2}\big)}
$$
and with $\Psi:\K\to \K$ defined by $\Psi a=\gamma_1^{1-1/p} \gamma_2^{1-1/u}a$ we have by \eqref{CC6} 
$$
\Psi I^{N_1,N_2}R=I^{M_1,M_2}.
$$
Now Proposition \ref{pro:2} and \eqref{CC8} yield that for all $n\in\N_0$ 
\begin{equation*}
e_n^\set(I^{M_1,M_2},B_{L_p^{M_1}(L_u^{M_2})})\le \gamma_1^{1-1/p} \gamma_2^{1-1/u}e_n^\set(I^{N_1,N_2},B_{L_p^{N_1}(L_u^{N_2})})\le 4e_n^\set(I^{N_1,N_2},B_{L_p^{N_1}(L_u^{N_2})}).
\end{equation*}
\end{proof}
\begin{corollary}
Let $1\le p,u\le \infty$ and let $\bar{p},\bar{u}$ be given by \eqref{CK5} and let $0<c_0<1$ be the constant from Proposition \ref{pro:1}. There exist constants $c_{1-5}>0$ such that for each $n,N_1,N_2\in \N$, with
\begin{equation}
\label{CF4}
n<\frac{c_0}{2}N_1N_2
\end{equation}
we have
\begin{eqnarray}
e_n^\ran(I^{N_1,N_2},B_{L_p^{N_1}(L_u^{N_2})})
&\ge&  c_1N_2^{1/u-1/\bar{p}}n^{1/\bar{p}-1}+c_1n^{-1/2}\quad (n\ge N_2)
\label{CD5}\\
e_n^\ran(I^{N_1,N_2},B_{L_p^{N_1}(L_u^{N_2})})
&\ge&  c_2N_1^{1/\bar{p}-1/u}n^{1/u-1}+c_2n^{-1/2}\quad (n\ge N_1)
\label{CD4}\\
e_n^\ranno(I^{N_1,N_2},B_{L_p^{N_1}(L_u^{N_2})})
&\ge& c_3 N_1^{1/p-1/2}n^{-1/2}\quad (n\ge N_1)
\label{CD7}\\
e_n^\ran(I^{N_1,N_2},B_{L_p^{N_1}(L_u^{N_2})})
&\ge&  c_4n^{1/\bar{p}-1}\quad (n< N_1)
\label{CD8}\\
e_n^\ran(I^{N_1,N_2},B_{L_p^{N_1}(L_u^{N_2})})
&\ge&  c_5n^{1/\bar{u}-1}\quad (n< N_2).
\label{CD9}
\end{eqnarray}
\end{corollary}

\begin{proof}
For all lower bounds we use the relations between average case and randomized setting \eqref{CK3}--\eqref{CK4}.
First we show \eqref{CD5}. We put
\begin{equation}
\label{CA5}
M_1=\left\lfloor\frac{2n}{c_0N_2}\right\rfloor +1,
\end{equation}
consequently 
\begin{equation}
\label{CA6}
n<\frac{c_0}{2}M_1N_2,
\end{equation}
and, using \eqref{CF4}, 
\begin{equation}
M_1\le \frac{2n}{c_0N_2}+1<  N_1+1,
\label{CJ0}
\end{equation}
hence
\begin{equation}
\label{CA8}
M_1\le N_1.
\end{equation}
Now assume $n\ge N_2$, then the first relation of \eqref{CJ0} implies
\begin{equation}
\label{CD3}
 M_1\le \frac{
2n+c_0N_2}{c_0N_2}\le \frac{(c_0+2)n}{c_0N_2}.
\end{equation}
By \eqref{CA8} we can apply Lemma \ref{lem:1} with $M_2=N_2$, giving  
\begin{eqnarray}
\label{CJ6}
 e_n^\ran(I^{N_1,N_2},B_{L_p^{N_1}(L_u^{N_2})})
 &\ge& \frac14 e_n^\ran(I^{M_1,N_2},B_{L_p^{M_1}(L_u^{N_2})}).
\end{eqnarray}
Taking into account \eqref{CA6}, we obtain from \eqref{UL2}--\eqref{CB2} of Proposition \ref{pro:1} and \eqref{CD3}
\begin{eqnarray*}
\lefteqn{e_n^\ran(I^{M_1,N_2},B_{L_p^{M_1}(L_u^{N_2})})}
\notag\\
&\ge&  \frac12\max\Big( e_{2n}^\avg(I^{M_1,N_2},\mu_{M_1,N_2}^{(1)}),e_{2n}^\avg(I^{M_1,N_2},\mu_{M_1,N_2}^{(2)}),e_{2n}^\avg(I^{M_1,N_2},\mu_{M_1,N_2}^{(3)})\Big)
\notag\\
&\ge& cM_1^{1/p-1}N_2^{1/u-1}+c(M_1N_2)^{-1/2}+cM_1^{-1/2}N_2^{1/u-1}
\notag\\
&\ge& cM_1^{1/\bar{p}-1}N_2^{1/u-1}+c(M_1N_2)^{-1/2}\ge   cN_2^{1/u-1/\bar{p}}n^{1/\bar{p}-1}+cn^{-1/2},
\end{eqnarray*}
which together with    \eqref{CJ6} proves \eqref{CD5}.

Next we verify \eqref{CD4} and \eqref{CD7}, where we define (symmetrically to the above)
\begin{equation}
\label{CA5E}
M_2=\left\lfloor\frac{2n}{c_0N_1}\right\rfloor +1,
\end{equation}
which analogously implies
\begin{eqnarray}
n&<&\frac{c_0}{2}N_1M_2\label{CA6E}\\
M_2&\le&  N_2.\label{CA8E}
\end{eqnarray}
Suppose that $n\ge N_1$, then according to \eqref{CA5E}
\begin{eqnarray}
M_2 &\le&  \frac{(c_0+2)n}{c_0N_1}.\label{CD3E}
\end{eqnarray}
Based on \eqref{CA8E}, we apply Lemma \ref{lem:1} with $M_1=N_1$ and get
\begin{eqnarray}
 e_n^\ran(I^{N_1,N_2},B_{L_p^{N_1}(L_u^{N_2})})
 &\ge& \frac14 e_n^\ran(I^{N_1,M_2},B_{L_p^{N_1}(L_u^{M_2})}).
\label{CJ2} 
\end{eqnarray}
Moreover, by \eqref{CA6E}, \eqref{UL2}--\eqref{CB2},  and \eqref{CD3E},
\begin{eqnarray*}
\lefteqn{ e_n^\ran(I^{N_1,M_2},B_{L_p^{N_1}(L_u^{M_2})})}
\notag\\
 &\ge& \frac12 \max\Big(e_{2n}^\avg(I^{N_1,M_2},\mu_{N_1,M_2}^{(1)}),e_{2n}^\avg(I^{N_1,M_2},\mu_{N_1,M_2}^{(2)}),e_{2n}^\avg(I^{N_1,M_2},\mu_{N_1,M_2}^{(3)})\big)
\notag\\
&\ge& cN_1^{1/\bar{p}-1}M_2^{1/u-1}+c(N_1M_2)^{-1/2}\ge cN_1^{1/\bar{p}-1/u}n^{1/u-1}+cn^{-1/2},
\end{eqnarray*}
which together with \eqref{CJ2} shows \eqref{CD4}. 
Furthermore, with Lemma \ref{lem:1} and \eqref{UQ3},
\begin{eqnarray*}
 e_n^\ranno(I^{N_1,N_2},B_{L_p^{N_1}(L_u^{N_2})}) &\ge& \frac14 e_n^\ranno(I^{N_1,M_2},B_{L_p^{N_1}(L_u^{M_2})})\ge\frac18 e_{2n}^\avgno(I^{N_1,M_2},\mu_{N_1,M_2}^{(4)})
\\
&\ge&  cN_1^{1/p-1}M_2^{-1/2}\ge cN_1^{1/p-1/2}n^{-1/2},
\end{eqnarray*}
establishing \eqref{CD7}.

If $n<N_1$, we let $M_2$ be defined by \eqref{CA5E}, thus $1\le M_2\le \left\lfloor\frac{2}{c_0}\right\rfloor +1$. Moreover, \eqref{CA6E} and \eqref{CA8E} are fulfilled. Hence by Lemma \ref{lem:1}  
\begin{eqnarray*}
e_n^\ran(I^{N_1,N_2},B_{L_p^{N_1}(L_u^{N_2})})
&\ge& \frac14 e_n^\ran(I^{N_1,M_2},B_{L_p^{N_1}(L_u^{M_2})}). 
\end{eqnarray*}
Now \eqref{CA8E} 
 together with the already shown relation \eqref{CD5} gives for $n\ge M_2$
\begin{eqnarray*}
 e_n^\ran(I^{N_1,M_2},B_{L_p^{N_1}(L_u^{M_2})}) &\ge& cM_2^{1/u-1/\bar{p}}n^{1/\bar{p}-1}\ge cn^{1/\bar{p}-1},
\end{eqnarray*}
consequently, 
\begin{eqnarray*}
e_n^\ran(I^{N_1,N_2},B_{L_p^{N_1}(L_u^{N_2})})
&\ge& cn^{1/\bar{p}-1}, 
\end{eqnarray*}
while for $n< M_2$ the result follows by monotonicity of the $e_n^\ran$ from $n=M_2$, 
proving \eqref{CD8}.

If $n<N_2$, we let $M_1$ be defined by \eqref{CA5}, hence $1\le M_1\le \left\lfloor\frac{2}{c_0}\right\rfloor +1$ and \eqref{CA6} and \eqref{CA8} are satisfied. 
Lemma \ref{lem:1} and   \eqref{CD4} imply for $n\ge M_1$ 
\begin{eqnarray*}
e_n^\ran(I^{N_1,N_2},B_{L_p^{N_1}(L_u^{N_2})})
&\ge& \frac14 e_n^\ran(I^{M_1,N_2},B_{L_p^{M_1}(L_u^{N_2})}) \ge  c_2M_1^{1/\bar{p}-1/u}n^{1/u-1}+cn^{-1/2}\ge cn^{1/\bar{u}-1}.
\end{eqnarray*}
Again, for $n< M_1$ the results follows by monotonicity, verifying \eqref{CD9}.

\end{proof}
\begin{theorem}
\label{theo:1}
Let $1 \le p, u\le \infty $ and let $\bar{p},\bar{u}$ be given by \eqref{CK5}. Then there exists constants $0<c_0<1$, $c_{1-13}>0$,  such that for $n,N_1,N_2\in\N$  with 
$n < c_0N_1N_2$ the following hold: \\
If $p\ge u$, then
\begin{eqnarray}
\label{P7}
c_1\min(N_2^{1/\bar{u}-1/\bar{p}}n^{-1+1/\bar{p}},n^{-1+1/\bar{u}})
&\le& 
e_n^\ran(I^{N_1,N_2},B_{L_p^{N_1}\big(L_u^{N_2}\big)})
\le e_n^{\rm ran-non}(I^{N_1,N_2},B_{L_p^{N_1}\big(L_u^{N_2}\big)}) 
\nonumber\\
&\le& c_2\min(N_2^{1/\bar{u}-1/\bar{p}}n^{-1+1/\bar{p}},n^{-1+1/\bar{u}}).
\end{eqnarray}
If $p<u\le 2$, then  
\begin{eqnarray}
\label{P8}
c_3\min(N_1^{1/p-1/u}n^{-1+1/u},n^{-1+1/p})
&\le& 
e_n^\ran(I^{N_1,N_2},B_{L_p^{N_1}\big(L_u^{N_2}\big)})
\le e_n^{\rm ran-non}(I^{N_1,N_2},B_{L_p^{N_1}\big(L_u^{N_2}\big)}) 
\nonumber\\
&\le& c_4 \min(N_1^{1/p-1/u}n^{-1+1/u},n^{-1+1/p}).
\end{eqnarray}
If $2\le p<u$, then  
\begin{eqnarray}
\label{A2}
c_5n^{-1/2}
&\le& 
e_n^\ran(I^{N_1,N_2},B_{L_p^{N_1}\big(L_u^{N_2}\big)})
\le e_n^{\rm ran-non}(I^{N_1,N_2},B_{L_p^{N_1}\big(L_u^{N_2}\big)}) 
\le c_6n^{-1/2}.
\end{eqnarray}
If $p<2<u$ and $n\ge N_1$, then  
\begin{eqnarray}
\label{A3}
&& c_7N_1^{1/p-1/u}n^{-1+1/u}+c_7n^{-1/2}  
\le e_n^\ran(I^{N_1,N_2},B_{L_p^{N_1}\big(L_u^{N_2}\big)}) 
\nonumber\\[.2cm]
&\le &
 c_8N_1^{1/p-1/u}n^{-1+1/u}\left(\log(N_1+1)\right)^{1-1/u}+c_8n^{-1/2}\left(\log(N_1+1)\right)^{1/2}
\end{eqnarray}
and
\begin{eqnarray}
c_9N_1^{1/p-1/2}n^{-1/2}
  &\le& e_n^{\rm ran-non}(I^{N_1,N_2},B_{L_p^{N_1}\big(L_u^{N_2}\big)})
\le c_{10} N_1^{1/p-1/2}n^{-1/2}.
\label{A4}
\end{eqnarray}
If $p<2<u$ and $n< N_1$, then   
\begin{eqnarray}
\label{A7}
c_{11}n^{-1+1/p}
&\le& 
e_n^\ran(I^{N_1,N_2},B_{L_p^{N_1}\big(L_u^{N_2}\big)})
\le e_n^{\rm ran-non}(I^{N_1,N_2},B_{L_p^{N_1}\big(L_u^{N_2}\big)}) 
\le c_{12}n^{-1+1/p}.
\end{eqnarray}
In the deterministic setting we have for $1\le p,u\le \infty$
\begin{equation}
\label{AF3}
c_{13}
\le
e_n^\de\big(I^{N_1,N_2},B_{L_p^{N_1}\big(L_u^{N_2}\big)}\big)
\le e_n^{\deno}\big(I^{N_1,N_2},B_{L_p^{N_1}\big(L_u^{N_2}\big)}\big) 
\le 1.
\end{equation}
\end{theorem}
\begin{proof}
The upper bound of \eqref{AF3} is a trivial consequence of \eqref{AF2}. All other upper bounds except that of \eqref{A3} follow from \eqref{AF6} and Proposition \ref{pro:4}, since the involved algorithm $A^{(2)}_n$ is non-adaptive.   
To prove \eqref{A3} we assume $n\ge N_1$. First we suppose that
$$
n<6N_1\left\lceil c(1)\log(N_1+1)\right\rceil,
$$ 
where $c(1)$ stands for the constant $c_1$ from Proposition \ref{pro:5}. Then \eqref{A4} implies 
\begin{eqnarray*}
\lefteqn{e_n^\ran(I^{N_1,N_2},B_{L_p^{N_1}\big(L_u^{N_2}\big)})
\le e_n^{\rm ran-non}(I^{N_1,N_2},B_{L_p^{N_1}\big(L_u^{N_2}\big)})}
\\
&\le&cN_1^{1/p-1/u}n^{-1+1/u}\left(\frac{n}{N_1}\right)^{1/2-1/u}\le cN_1^{1/p-1/u}n^{-1+1/u} \left(\log(N_1+1)\right)^{1/2-1/u},
\end{eqnarray*}
which gives the upper bound of \eqref{A3}. Now assume
\begin{equation}
\label{C11}
n\ge 6N_1\left\lceil c(1)\log(N_1+1)\right\rceil.
\end{equation}
We set 
\begin{equation*}
m=\left\lceil c(1)\log(N_1+1)\right\rceil, \quad\tilde{n}=\left\lfloor\frac{n}{6\left\lceil c(1)\log(N_1+1)\right\rceil}\right\rfloor,
\end{equation*}
which implies $\tilde{n}\ge  N_1$, so that we can use Proposition \ref{pro:5} with $\tilde{n} $ instead of $n$. Hence by \eqref{WM1}
\begin{equation*}
\ca(A_{\tilde{n},\omega}^3)\le 6m\tilde{n}\le n,
\end{equation*}
and therefore
\begin{eqnarray}
\label{C2}
e_n^\ran(I^{N_1,N_2},B_{L_p^{N_1}\big(L_u^{N_2}\big)}) 
&\le & 
 c\left(N_1^{1/p-1/u}\tilde{n}^{-1+1/u}+\tilde{n}^{-1/2}\right).
\end{eqnarray}
From \eqref{C11} we obtain, using $\lfloor a\rfloor\ge a/2$ for $a\ge 1$,
\begin{eqnarray*}
\tilde{n}
&>&\frac{n}{12\left\lceil c(1)\log(N_1+1)\right\rceil}
\ge \frac{n}{12(c(1)+1)\log(N_1+1)},
\end{eqnarray*}
which combined with \eqref{C2} completes the proof of  the upper bound in \eqref{A3}.

Now we prove the lower bounds. We start with the case $p\ge u$, thus with  \eqref{P7}. Relation  \eqref{CD9} gives the lower bound for $n<N_2$. The case $n\ge N_2$ follows from \eqref{CD5}. Indeed, if $u\ge 2$, we have $\bar{p}=\bar{u}=2$ and the second term on the right-hand side of \eqref{CD5} gives the lower bound of \eqref{P7}, while for $u<2$ the first term yields the desired result.   

Next suppose that $p<u$. The lower bound of \eqref{P8} follows from \eqref{CD4} and \eqref{CD8}, while that of \eqref{A2} is a consequence of \eqref{CD5} and \eqref{CD9} (since $1/\bar{u}-1\ge -1/2$ for all $u$). Furthermore,  \eqref{CD4} implies  
the lower bound in \eqref{A3}, while  \eqref{CD7} implies that of \eqref{A4}, and \eqref{CD8} yields that of \eqref{A7}.

Finally, in the deterministic setting we use 
$
B_{L_\infty^{N_1}\big(L_\infty^{N_2}\big)}\subseteq B_{L_p^{N_1}\big(L_u^{N_2}\big)},
$
hence 
\begin{equation}
\label{AF4}
e_n^\de\big(I^{N_1,N_2},B_{L_\infty^{N_1}\big(L_\infty^{N_2}\big)}\big)\le e_n^\de\big(I^{N_1,N_2},B_{L_p^{N_1}\big(L_u^{N_2}\big)}\big).
\end{equation}
On the other hand, it is well-known that
\begin{equation*}
e_n^\de\big(I^{N_1,N_2},B_{L_\infty^{N_1}\big(L_\infty^{N_2}\big)}\big)=e_n^\de\big(I^{N_1N_2},B_{L_\infty^{N_1N_2}}\big) \ge c,
\end{equation*}
which together with \eqref{AF4} and \eqref{AF5} gives \eqref{AF3}.

\end{proof}
We want to estimate the order of the largest gap between non-adaptive and adaptive minimal errors among the family of studied integration problems. For this purpose we carry out an analysis similar to that in \cite{Hei23a}, see Corollary 4.6 there.   Let $p<2<u$.
Consider the region $N_1\le n<c(0)N_1N_2$, with $0<c(0)<1$ standing for the constant $c_0$ from Theorem \ref{theo:1} and define 
\begin{eqnarray*}
\gamma(p,u,n,N_1,N_2)&=& \frac{e_n^\ranno(I^{N_1,N_2},B_{L_p^{N_1}(L_u^{N_2})})}{e_n^\ran(I^{N_1,N_2},B_{L_p^{N_1}(L_u^{N_2})})}\\[.2cm]
\gamma(p,u,n)&=&\sup_{N_1,N_2:\,N_1\le n<c(0)N_1N_2}\gamma(p,u,n,N_1,N_2).
\end{eqnarray*}
\begin{corollary}
\label{cor:3}
Let $1\le p<2<u\le \infty$. Then there are constants $c_1,c_2>0$ such that for all $n\in\N$
\begin{equation}
\label{N1}
c_1(\log(n+1))^{1/u-1}n^{\frac{\left(\frac{1}{p}-\frac{1}{2}\right)\left(\frac{1}{2}-\frac{1}{u}\right)}{\frac{1}{p}-\frac{1}{u}}}\le \gamma(p,u,n)\le c_2n^{\frac{\left(\frac{1}{p}-\frac{1}{2}\right)\left(\frac{1}{2}-\frac{1}{u}\right)}{\frac{1}{p}-\frac{1}{u}}}.
\end{equation}
The exponent of $n$ attains its maximal value $1/4$ iff $p=1$, $u=\infty$. In this case 
the following holds. For any constants $c_3,c_4$  with $c(0)^{1/2}<c_3<c_4$ there are constant $c_5,c_6>0$ such that for all $n\in\N$ with $n\ge c_4^2$  and all $N_1,N_2\in [c_3 n^{1/2},c_4n^{1/2}]$ \; 
\begin{equation}
 c_5(\log(n+1))^{-1}n^{1/4}\le\frac{e_n^{\rm ran-non }(I^{N_1,N_2},B_{L_1^{N_1}(L_\infty^{N_2})})}{e_n^{\rm ran }(I^{N_1,N_2},B_{L_1^{N_1}(L_\infty^{N_2})})}\le c_6n^{1/4}.
\label{CH9}
\end{equation}
\end{corollary}
\begin{proof}
We will estimate $\gamma(p,u,n)^{-1}$. By \eqref{A3} and \eqref{A4} of  Theorem \ref{theo:1} there are constants $c_1,c_2>0$ such that for $n,N_1,N_2\in\N$ with $N_1\le n<c(0)N_1N_2$
\begin{eqnarray}
&& c_1 \max\bigg(\Big(\frac{n}{N_1}\Big)^{1/u-1/2}, 
N_1^{1/2-1/p} \bigg)
\le\gamma(p,q,n,N_1,N_2)^{-1}
\notag\\
&\le&c_2\Big(\frac{n}{N_1}\Big)^{1/u-1/2}(\log(N_1+1))^{1-1/u}
+c_2 N_1^{1/2-1/p} (\log(N_1+1))^{1/2}
\label{CI1}
\end{eqnarray}
and therefore
\begin{eqnarray}
&& c_1 \min_{N_1: 1\le N_1\le n}\max\bigg(\Big(\frac{n}{N_1}\Big)^{1/u-1/2}, 
N_1^{1/2-1/p} \bigg)
\le\gamma(p,q,n)^{-1}
\notag\\
&\le&c_2\min_{N_1,N_2:\,N_1\le n<c(0)N_1N_2}\bigg(\Big(\frac{n}{N_1}\Big)^{1/u-1/2}(\log(N_1+1))^{1-1/u}
\notag\\
&&\hspace{4.5cm}+ N_1^{1/2-1/p} (\log(N_1+1))^{1/2}\bigg).
\label{M5}
\end{eqnarray}
Define $x_0$ by
\begin{equation}
\Big(\frac{n}{x_0}\Big)^{1/u-1/2}= 
x_0^{1/2-1/p},
\label{N8}
\end{equation}
then
\begin{equation}
x_0=n^{\frac{\frac{1}{u}-\frac{1}{2}}{\frac{1}{u}-\frac{1}{p}}}, \quad x_0\in[1,n],\label{N7}
\end{equation}
and 
\begin{equation*}
\min_{x\in[1,n]}\max\left(\left(\frac{n}{x}\right)^{1/u-1/2}, x^{1/2-1/p}\right)= x_0^{1/2-1/p} =n^{-\frac{\left(\frac{1}{p}-\frac{1}{2}\right)\left(\frac{1}{2}-\frac{1}{u}\right)}{\frac{1}{p}-\frac{1}{u}}},
\end{equation*}
which together with the lower bound in \eqref{M5} shows 
the upper bound in \eqref{N1}. 

Exactly as in \cite{Hei23a} we define
\begin{equation*}
N_1=\left\lceil x_0\right\rceil,\quad 
N_2=\left\lfloor \frac{n}{c(0)x_0}\right\rfloor+1,
\end{equation*}
and conclude 
\begin{equation*}
N_1\le n,\quad x_0\le N_1<2x_0,\quad \frac{n}{c(0)N_1}\le \frac{n}{c(0)x_0}<N_2<\frac{2n}{c(0)x_0}\le \frac{2n}{c(0)}\,,
\end{equation*}
in particular we have $N_1\le n<c(0)N_1N_2$. Consequently, the upper bound of  \eqref{M5} together with \eqref{N8} and  \eqref{N7} gives 
\begin{eqnarray*}
 &&\gamma(p,q,n)^{-1}
\\
&\le&c\left(\left(\frac{n}{2x_0}\right)^{1/u-1/2}(\log(2x_0+2c(0)^{-1}n)^{1-1/u}+ 
x_0^{1/2-1/p}(\log(2x_0+2c(0)^{-1}n)^{1/2} \right)
\\
&\le&c n^{-\frac{\left(\frac{1}{p}-\frac{1}{2}\right)\left(\frac{1}{2}-\frac{1}{u}\right)}{\frac{1}{p}-\frac{1}{u}}}
(\log(n+1))^{1-1/u},
\end{eqnarray*}
implying the lower bound of \eqref{N1}.

The exponent of the  gap between non-adaption and adaption satisfies
\begin{eqnarray*}
\frac{\left(\frac{1}{p}-\frac{1}{2}\right)\left(\frac{1}{2}-\frac{1}{u}\right)}{\frac{1}{p}-\frac{1}{u}}
&\le&
\frac{\left(\frac{1}{p}-\frac{1}{2}+\frac{1}{2}-\frac{1}{u}\right)^2}{4\left(\frac{1}{p}-\frac{1}{u}\right)}
=\frac{\frac{1}{p}-\frac{1}{u}}{4}\le \frac{1}{4}.
\end{eqnarray*}
It follows that the exponent attains the maximal value 1/4 iff $p=1$, $u=\infty$. With this choice and $N_1(n), N_2(n)$ as assumed, relation \eqref{CH9} is a direct consequence of \eqref{CI1}.

\end{proof}
\section{Passing from finite to infinite dimensional problems}
\label{sec:5}
So far we produced integration problems with suitable gaps between adaptive and non-adaptive randomized $n$-th minimal errors for each $n$ separately. This raises the question if there are infinite dimensional examples with respective gaps for all $n$ simultaneously.
In this section we use infinite direct sums of finite dimensional spaces to produce such examples. We study various aspects of the adaption problem in these spaces including operators represented as an infinite sum of finite dimensional operators and similarly, operators with a respective diagonal representation.

Let $X_k$ $(k\in\N_0)$ be Banach spaces, $1\le p_1\le \infty$, and let 
\begin{equation}
\label{CF7}
X=\bigg(\bigoplus_{k=0}^\infty X_k\bigg)_{p_1}
\end{equation}
be the space of tuples $x=(x_k)_{k=0}^\infty$ with $x_k\in X_k$, endowed with the norm  
$$
\|x\|=\big\|(\|x_k\|_{X_k})_{k=0}^\infty\big\|_{\ell_{p_1}}=\left\{\begin{array}{lll}
   & \displaystyle\bigg(\sum_{k=0}^\infty \|x_k\|_{X_k}^{p_1}\bigg)^{1/{p_1}}\quad\mbox{if}\quad  {p_1}<\infty  \\[.4cm]
   & \displaystyle\sup_{0\le k<\infty}\|x_k\|_{X_k}\quad\mbox{if}\quad  {p_1}=\infty.    
    \end{array}
\right.
$$ 
We define the following mappings
\begin{eqnarray}
&& J_k:X_k\to X, \quad J_kx_k=(\delta^{kl}x_k)_{l=0}^\infty 
\label{PD1}\\
&&Q_k:X\to X_k, \quad Q_k((x_l)_{l=0}^\infty)=x_k  
\label{PC7}\\
&&P_k:X\to X, \quad P_k=J_kQ_k,
\label{PC8}
\end{eqnarray}
where $\delta^{kl}$ stands for the Kronecker symbol.  Furthermore, let $\emptyset\ne \Lambda_k\subseteq X_k^*$ and 
\begin{eqnarray}
\Lambda=\bigcup_{k=0}^\infty \Phi_k(\Lambda_k),\text{ where }\Phi_k:\Lambda_k\to X^*, \quad (\Phi_k(\lambda_k))(x)=\lambda_k(x_k).\label{CF3}
\end{eqnarray}

\subsection{Infinite sum problems}
\label{sec:5.1}
Let $G$ be a Banach space and  consider the problems $\mathcal{P}_k=(B_{X_k},G,S_k,\K,\Lambda_k)$, where $X_k,\Lambda_k$ are as above, while $S_k\in L(X_k,G)$ $(k\in\N_0)$. We 
define the sum problem $\mathcal{P}=(B_X,G,S,\K,\Lambda)$ by
\begin{eqnarray}
&& S:X\to G,\quad S((x_k)_{k=0}^\infty)=\sum_{k=0}^\infty S_kx_k,\label{MA3}
\end{eqnarray}
with  $X$ and $\Lambda$ as given in \eqref{CF7} and \eqref{CF3}. A sufficient condition for $S$ being well-defined and bounded is 
\begin{equation}
\label{PC4}
\|(\|S_k\|)\|_{\ell_{p_1^*}}<\infty,
\end{equation}
where $1/p_1^*=1-1/{p_1}$.
Indeed, 
\begin{eqnarray}
\|S\|&\le&\sup_{\|(\|x_k\|)\|_{\ell_{p_1}}\le 1}\sum_{k=0}^\infty \|S_k(x_k)\|
=\sup_{\|(b_k)\|_{\ell_{p_1}}\le 1}\sum_{k=0}^\infty \sup_{\|x_k\|\le |b_k|}\|S_k(x_k)\|
\nonumber\\
&=&\sup_{\|(b_k)\|_{\ell_{p_1}}\le 1}\sum_{k=0}^\infty \|S_k\||b_k|=\|(\|S_k\|)\|_{\ell_{p_1^*}}<\infty.\label{AH0}
\end{eqnarray}
In a similar way it follows from \eqref{PC4} that
\begin{eqnarray}
\label{PC5}
\bigg\|S-\sum_{j=0}^{k_1}SP_k\bigg\|\le\|( \|S_k\|)_{k=k_1+1}^\infty\big\|_{\ell_{p_1^*}}\quad(k_1\in\N_0).
\end{eqnarray}
Moreover, if $G=\K$, then equality holds in the first relation of \eqref{AH0} and in \eqref{PC5}.
\begin{lemma}
\label{lem:6}
Let ${\rm set}\in\{\de,\deno,\ran,\ranno\}$ and assume that \eqref{PC4} holds. Then for $n\in\N_0$
\begin{eqnarray}
\label{PA2}
 e_n^\set(S_k,B_{X_k}, G)\le e_n^\set(S,B_X,G).
\end{eqnarray}
Moreover, for $k_1\in\N_0$, $(n_k)_{k=0}^{k_1}\subset\N_0$, and $n$ defined by $n=\sum_{k=0}^{k_1}n_k$
\begin{eqnarray}
e_n^\set(S,B_X,G)
&\le& \big\|S-S\sum_{k=0}^{k_1} P_k\big\|_{L(X,G)}+\sum_{k=0}^{k_1} e_{n_k}^\set(S_k,B_{X_k}, G)
\label{PA3}\\
&\le& \big\|( \|S_k\|)_{k=k_1+1}^\infty\big\|_{\ell_{p_1^*}}+\sum_{k=0}^{k_1} e_{n_k}^\set(S_k,B_{X_k}, G).
\label{PD3}
\end{eqnarray}
\end{lemma}
\begin{proof} 
To prove \eqref{PA2} we  show that $\mathcal{P}_k$ reduces to $\mathcal{P}$. We have 
\begin{equation}
\label{PD2}
SJ_k=S_k,\quad J_k(B_{X_k})\subseteq B_X.
\end{equation}
Fix any $\lambda_{k,0}\in\Lambda_k$ and  define $\eta_1:\Lambda\to\Lambda_k$, $\rho:\Lambda\times \K \to \K$ by setting for $\lambda\in\Lambda$ and $a\in \K$
\begin{eqnarray*}
\eta_1(\lambda)&=&\left\{\begin{array}{lll}
   \lambda_k & \quad\mbox{if}\quad \lambda=\Phi_k(\lambda_k) \text{ for some } \lambda_k\in\Lambda_k  \\
  \lambda_{k,0} & \quad\mbox{otherwise,} 
    \end{array}
\right. 
\\
\rho(\lambda,a)&=&\left\{\begin{array}{lll}
   a & \quad\mbox{if}\quad \lambda\in\Phi_k(\Lambda_k)   \\
  0 & \quad\mbox{if}\quad \lambda\not\in\Phi_k(\Lambda_k).  
    \end{array}
\right. 
\end{eqnarray*}
We show that for all $\lambda\in\Lambda$ and all $x_k\in B_{X_k}$
\begin{equation}
\label{PB6a}
\lambda(J_k(x_k))=\rho(\lambda, (\eta_1(\lambda))(x_k)).
\end{equation}
If $\lambda=\Phi_k(\lambda_k)$  for some $\lambda_k\in\Phi_k(\Lambda_k) $, then 
\begin{equation*}
\lambda(J_k(x_k))=\lambda_k(x_k)=\rho(\lambda, \lambda_k(x_k))=\rho(\lambda, (\eta_1(\lambda))(x_k)).
\end{equation*}
On the other hand, if $\lambda\not\in\Phi_k(\Lambda_k)$, thus $\lambda=\Phi_l(\lambda_l)$  for some  $\lambda_l\in\Lambda_l$ with $l\ne k$, we have
\begin{equation*}
\lambda(J_k(x_k))=\Phi_l(\lambda_l)(J_k(x_k))=0=\rho(\lambda, \lambda_{k,0}(x_k))=\rho(\lambda, (\eta_1(\lambda))(x_k)).
\end{equation*}
This shows \eqref{PB6a}. Based on \eqref{PD2} and \eqref{PB6a} we conclude from Proposition  \ref{pro:2}
\begin{eqnarray*}
e_n^\set(S_k,B_{X_k},G)&\le& e_n^\set(S,B_X,G).
\end{eqnarray*}

Next we turn to \eqref{PA3} and show  that  $(B_X,G,SP_k,\K,\Lambda)$ reduces to $\mathcal{P}_k$. 
Here we have 
\begin{equation}
\label{PB9}
SP_k= S_k Q_k, \quad Q_k(B_X)= B_{X_k}.
\end{equation}
Set 
$$
\eta_1=\Phi_k:\Lambda_k\to \Lambda,\quad \rho:\Lambda_k\times \K\to \K,\quad
\rho(\lambda_k,a)=a.
$$ 
It follows that for $\lambda_k\in\Lambda_k$ and $x=(x_k)\in X$
\begin{eqnarray}
\label{PB8}
\lambda_k(Q_kx)&=&\lambda_k(x_k)=\rho(\lambda_k,\lambda_k(x_k))
=\rho(\lambda_k,(\Phi_k(\lambda_k))(x))
\nonumber\\
&=&\rho(\lambda_k, (\eta_1(\lambda_k))(x)).
\end{eqnarray}
Using Proposition \ref{pro:2} again we obtain from \eqref{PB9} and \eqref{PB8}
\begin{eqnarray}
 e_n^\set(SP_k,B_X,G)&\le&e_n^\set(S_k,B_{X_k},G).\label{CF5}
\end{eqnarray}
Consequently, now using Proposition \ref{pro:3} and \eqref{CF5},
\begin{eqnarray*}
\label{PC6}
e_n^\set(S,B_X,G)&\le& \big\|S-S\sum_{k=0}^{k_1} P_k\big\|_{L(X,G)}+e_n^\set\bigg(\sum_{k=0}^{k_1}SP_k,B_X, G\bigg)
\\
&\le& \big\|S-S\sum_{k=0}^{k_1} P_k\big\|_{L(X,G)}+\sum_{k=0}^{k_1}e_{n_k}^\set(SP_k,B_X, G)
\\
&\le& \big\|S-S\sum_{k=0}^{k_1} P_k\big\|_{L(X,G)}+\sum_{k=0}^{k_1} e_{n_k}^\set(S_k,B_{X_k}, G_k),
\end{eqnarray*}
which is \eqref{PA3}, and combined with \eqref{PC5} gives \eqref{PD3}.

\end{proof}

We apply Lemma \ref{lem:6} to the following situation. Let $1\le p,u\le \infty$, $\alpha\in\R$, $\alpha>1$,
and for $k\in\N_0$ 
\begin{equation*}
N_k=2^k,\quad X_k=L_p^{N_k}(L_u^{N_k}),\quad G=\K,
\quad S_k=2^{-\alpha k}I^{N_k,N_k},
\end{equation*}
and $\mathcal{P}_k=\big(B_{L_p^{N_k}(L_u^{N_k})},\K,2^{-\alpha k}I^{N_k,N_k},\K,\Lambda^{N_k,N_k}\big) $, where $\Lambda^{N_k,N_k}$ is standard information, see \eqref{CI7}. Thus we consider the problem defined by \eqref{CF7}, \eqref{CF3}, and \eqref{MA3}
\begin{equation*}
X=\bigg(\bigoplus_{k=0}^\infty L_p^{N_k}(L_u^{N_k})\bigg)_{\ell_{p_1}}, \quad G=\K,\quad \Lambda=\bigcup_{k=0}^\infty \Phi_k(\Lambda^{N_k,N_k}),
\end{equation*}
\begin{equation*}
S=I:X\to \K,\quad I(f_k)_{k=0}^\infty=\sum_{k=0}^\infty 2^{-\alpha k}I^{N_k,N_k}f_k
=\sum_{k=0}^\infty 2^{-(\alpha+2) k}\sum_{i,j=1}^{N_k}f_k(i,j)\in \K.
\end{equation*}
Since $\|I^{N_k,N_k}:L_p^{N_k}(L_u^{N_k})\to \K\|=1$, we get from \eqref{AF2}, \eqref{PC4}, and the assumption $\alpha>1$,  
\begin{eqnarray*}
\|I:X\to\K\| &=&\|(2^{-\alpha k}\|I^{N_k,N_k}\|)\|_{\ell_{p_1^*}}<\infty
\end{eqnarray*}
and from \eqref{PC5} for $k_1\in\N_0$
\begin{equation}
\label{CK1}
\bigg\|I-\sum_{j=0}^{k_1}IP_k\bigg\|=\big\|(2^{-\alpha k} \|I^{N_k,N_k}\|)_{k=k_1+1}^\infty\big\|_{\ell_{p_1^*}}\le c2^{-\alpha k_1}
\end{equation}
for some constant $c>0$. Now we transfer Theorem \ref{theo:1} to this infinite dimensional situation. We only consider the case $1\le p<2<u\le\infty$, the other cases can be treated analogously.
\begin{proposition}
\label{pro:6} Let $1\le p<2<u\le \infty$ and $\alpha>1$. Then there are constants $c_{1-6}>0$ and $n_0\in\N$ such that for $n\ge n_0$
\begin{eqnarray*}
c_1n^{-\alpha/2+1/(2p)+1/(2u)-1}+c_1n^{-\alpha/2-1/2}&\le& e_n^\ran(I,B_X)
\le c_2n^{-\alpha/2+1/(2p)+1/(2u)-1}(\log(n+1))^{1-1/u}
\notag\\[.2cm]
&&+c_2n^{-\alpha/2-1/2}(\log(n+1))^{1/2}
 \\[.3cm]
c_3n^{-\alpha/2+1/(2p)-3/4}&\le& e_n^\ranno(I,B_X)
\le c_4n^{-\alpha/2+1/(2p)-3/4}
\\[.3cm]
c_5n^{-\alpha/2}&\le& e_n^\de(I,B_X)\le e_n^\deno(I,B_X)
\le c_6n^{-\alpha/2}.
\end{eqnarray*}
\end{proposition}
\begin{proof}
Let $0<c(0)<1$ be the constant $c_0$ from Theorem \ref{theo:1}. First we show the lower bounds.  Set 
\begin{equation}
\label{PB5}
c(1)=\lceil -\log_2 c(0)\rceil+1\ge 2,
\end{equation}
thus $c(0)\ge 2^{-c(1)+1}$.
For any $k\in\N$, with 
\begin{equation}
\label{PB1}
k\ge c(1),\quad \text{thus }\quad c(0)2^k\ge 2^{k-c(1)+1}\ge 2,
\end{equation}
we set $m_k :=\lceil  c(0)2^{2k}\rceil-1$, hence
\begin{equation}
\label{PB4}
m_k <c(0)2^{2k}=c(0)N_k^2, \quad  m_k \ge c(0)2^{2k-1}\ge 2^{2k-c(1)}\ge 2^{k}=N_k.
\end{equation}
From  \eqref{A3}, \eqref{A4}, and \eqref{AF3} of Theorem \ref{theo:1}, \eqref{PA2}, and \eqref{PB4} we obtain 
\begin{eqnarray}
e_{m_k}^\ran(I,B_X)&\ge& 2^{-\alpha k} e_{m_k}^{\rm ran }(I^{N_k,N_k},B_{L_p^{N_k}(L_u^{N_k})})
\notag\\
&\ge& c2^{-\alpha k}(N_k^{1/p-1/u}m_k^{-1+1/u}+m_k^{-1/2})
\ge  c2^{-\alpha k}\big(2^{(1/p+1/u-2)k}+2^{-k}\big)\label{AD9}
\\[.2cm]
e_{m_k}^\ranno(I,B_X)&\ge& 2^{-\alpha k} e_{m_k}^{\rm ran-non }(I^{N_k,N_k},B_{L_p^{N_k}(L_u^{N_k})})
\notag\\
&\ge& c2^{-\alpha k}N_k^{1/p-1/2} m_k^{-1/2}\ge c2^{(-\alpha +1/p-3/2)k}
\label{PA1}
\\[.2cm]
e_{m_k}^\de(I,B_X)&\ge& 2^{-\alpha k} e_{m_k}^{\de }(I^{N_k,N_k},B_{L_p^{N_k}(L_u^{N_k})})
\ge c2^{-\alpha k}.\label{AG0}
\end{eqnarray}
We transform  \eqref{AD9}--\eqref{AG0} into estimates for arbitrary $n\in \N$ with $n\ge 2^{c(1)}$. Choose a $k\ge c(1)$ such that 
\begin{equation}
2^{2k-c(1)}\le n<2^{2(k+1)-c(1)},
\label{CK2}
\end{equation}
then by \eqref{PB4} $n<m_{k+1}$ and  \eqref{AD9}--\eqref{CK2} imply for $n\ge 2^{c(1)}$ 
\begin{eqnarray}
e_n^\ran(I,B_X)&\ge&e_{m_k}^\ran(I,B_X)\ge cn^{-\alpha/2+1/(2p)+1/(2u)-1}+cn^{-\alpha/2-1/2}\quad 
\label{AE0}
\\
e_n^\ranno(I,B_X)&\ge& e_{m_k}^\ranno(I,B_X)\ge cn^{-\alpha/2+1/(2p)-3/4}
\label{CJ9}\\
e_n^\de(I,B_X)&\ge& e_{m_k}^\de(I,B_X)\ge c2^{-\alpha k}\ge cn^{-\alpha/2}. \label{AF9}
\end{eqnarray}

Next we turn to the upper estimates.  Let $\beta$ and $\delta$ be reals satisfying
\begin{eqnarray}
&&1<\beta<2, \label{AE3}
\\
&&0<\delta(\beta-1)<2-\beta,
\label{CI9}
\end{eqnarray}
and let $c(2)\in\N$ be a constant chosen sufficiently large, so  that 
\begin{equation*}
\lceil  c(0)2^{(2-\delta(\beta-1))k}\rceil-1\ge 2^{\beta k}  
\end{equation*}
for all $k\ge c(2)$.
Now fix $k_0\in\N$ with  $k_0\ge c(2)$ and   
define 
\begin{equation}
\label{PA0}
k_1=\lfloor \beta k_0\rfloor, \quad  n_k=\left\{\begin{array}{lll}
 2^{2k}  & \quad\mbox{if}\quad  k<k_0  \\
   \lceil  c(0)2^{2k_0-\delta(k-k_0)}\rceil-1& \quad\mbox{if}\quad k_0\le k\le k_1.
    \end{array}
\right. 
\end{equation}
For $k_0\le k\le k_1$ we have, taking into account \eqref{CI9},
\begin{eqnarray}
c(0)N_k^2=c(0)2^{2k}> n_k&\ge& n_{k_1}=  \lceil  c(0)2^{2k_0-\delta(k_1-k_0)}\rceil-1
\notag\\
&\ge& 
\lceil  c(0)2^{(2-\delta(\beta-1))k_0}\rceil-1\ge 2^{\beta k_0}\ge 2^{k_1}\ge N_k
\label{CK0}\\[.2cm]
n_k&\ge&  c(0)2^{2k_0-\delta(k-k_0)-1}, \label{CJ8}
\end{eqnarray}
the latter being a consequence of $ \lceil  a \rceil-1\ge a/2$ for $a>1$.
Relation \eqref{PA0} implies 
\begin{eqnarray}
\sum_{k=0}^{k_1}n_k&\le& c(3) 2^{2k_0}
\label{PA6}
\end{eqnarray}
for some constant $c(3)\in\N$
and 
\begin{eqnarray}
e_{n_k}^\ran(I^{N_k,N_k},B_{L_p^{N_k}(L_u^{N_k})})&=&e_{n_k}^\ranno(I^{N_k,N_k},B_{L_p^{N_k}(L_u^{N_k})})
\nonumber\\
&=&e_{n_k}^\de(I^{N_k,N_k},B_{L_p^{N_k}(L_u^{N_k})})=0\qquad(0\le k< k_0).
\label{PC0}
\end{eqnarray}
By \eqref{A3}, \eqref{A4}, and \eqref{AF3}   of Theorem \ref{theo:1},  and \eqref{CJ8}, for $k_0\le k\le k_1$ 
\begin{eqnarray}
e_{n_k}^\ran(I^{N_k,N_k},B_{L_p^{N_k}(L_u^{N_k})})&\le& c(k+1)^{1-1/u}N_k^{1/p-1/u}n_k^{-1+1/u}+c(k+1)^{1/2}n_k^{-1/2}
\notag\\[.1cm]
&\le&  c(k+1)^{1-1/u}2^{(1/p-1/u)k+(2/u-2)k_0+(1-1/u)\delta(k-k_0)}\notag\\
&&\hspace{2cm} +c(k+1)^{1/2}2^{-k_0+\delta(k-k_0)/2}
\label{PC1}
\\[.2cm]
e_{n_k}^\ranno(I^{N_k,N_k},B_{L_p^{N_k}(L_u^{N_k})})&\le& cN_k^{1/p-1/2} n_k^{-1/2}\le c2^{(1/p-1/2)k-k_0+\delta(k-k_0)/2}
\label{AE1}
\\[.2cm]
e_{n_k}^\deno(I^{N_k,N_k},B_{L_p^{N_k}(L_u^{N_k})})&\le& 1.\label{AG1}
\end{eqnarray}
Now we set
\begin{equation}
\beta=\frac{\alpha+1}{\alpha}.\label{CL0}
\end{equation}
Since $\alpha>1$,  \eqref{AE3} is satisfied and \eqref{CI9} turns into 
\begin{equation}
0<\delta<\alpha-1.\label{CK8}
\end{equation}
(The reason for not fixing this value of $\beta$ from the beginning is that we want to reuse relations \eqref{AE3}--\eqref{PA6} in the proof of Proposition
\ref{pro:7} with another value of $\beta$.) From \eqref{PD3}, \eqref{CK1}, \eqref{PA0}, and \eqref{PA6}--\eqref{CK8} we obtain
\begin{eqnarray}
\lefteqn{ e_{c(3) 2^{2k_0}}^\ran(I,B_X)}
\notag\\
&\le& 
c2^{-\alpha k_1} +c\sum_{k=k_0}^{k_1}2^{-\alpha k}\big((k+1)^{1-1/u}2^{(1/p-1/u)k+(2/u-2)k_0+(1-1/u)\delta(k-k_0)}
\nonumber\\
&&\hspace{4cm}+(k+1)^{1/2}2^{-k_0+\delta(k-k_0)/2}\big)
\nonumber\\
&\le& c2^{-\alpha \beta k_0} +c(k_0+1)^{1-1/u}2^{(-\alpha+1/p+1/u-2) k_0}\sum_{k=k_0}^{k_1}2^{(-\alpha+1/p-1/u+(1-1/u)\delta)(k-k_0)}
\nonumber\\
&&+c(k_0+1)^{1/2}2^{-(\alpha+1) k_0}\sum_{k=k_0}^{k_1}2^{(-\alpha+\delta/2)(k-k_0)}
\nonumber\\
&\le&c 2^{-(\alpha+1) k_0}+c(k_0+1)^{1-1/u}2^{(-\alpha+1/p+1/u-2) k_0}+c(k_0+1)^{1/2}2^{-(\alpha+1) k_0}
\nonumber\\[.2cm]
&\le&c(k_0+1)^{1-1/u}2^{(-\alpha+1/p+1/u-2) k_0}+c(k_0+1)^{1/2}2^{-(\alpha+1) k_0},
\label{PA5}
\end{eqnarray}
\begin{eqnarray}
e_{c(3) 2^{2k_0}}^\ranno(I,B_X)&\le& 
c2^{-\alpha k_1} +c\sum_{k=k_0}^{k_1}2^{-\alpha k +(1/p-1/2)k-k_0+\delta(k-k_0)/2}
\nonumber\\
&\le& 
c2^{-\alpha \beta k_0} +c2^{(-\alpha+1/p-3/2) k_0}\sum_{k=k_0}^{k_1}2^{(-\alpha+1/p-1/2+\delta/2)(k-k_0)}
\nonumber\\
&\le&c 2^{-(\alpha+1) k_0}+c2^{(-\alpha+1/p-3/2) k_0}\le c2^{(-\alpha+1/p-3/2) k_0},
\label{AE2}
\end{eqnarray}
\begin{eqnarray}
e_{c(3) 2^{2k_0}}^\deno(I,B_X)&\le& 
c2^{-\alpha k_1} +c\sum_{k=k_0}^{k_1}2^{-\alpha k}\le c2^{-\alpha k_0}.
\label{AG2}
\end{eqnarray}
Finally, let
\begin{equation*}
n\ge c(3)2^{2c(2)}.
\end{equation*}
Then we can choose $k_0\ge c(2)$ so that 
\begin{equation*}
c(3)2^{2k_0}\le n<c(3)2^{2(k_0+1)},
\end{equation*}
which together with \eqref{PA5}--\eqref{AG2}  gives for $n\ge c(3)2^{2c(2)}$
\begin{eqnarray*}
e_n^\ran(I,B_X)&\le& c(k_0+1)^{1-1/u}2^{(-\alpha+1/p+1/u-2) k_0}+c(k_0+1)^{1/2}2^{-(\alpha+1) k_0}
\notag\\[.15cm]
&\le& cn^{-\alpha/2+1/(2p)+1/(2u)-1}(\log(n+1))^{1-1/u}+cn^{-\alpha/2-1/2}(\log(n+1))^{1/2}
\\[.25cm]
e_n^\ranno(I,B_X)&\le& c2^{(-\alpha+1/p-3/2) k_0}\le cn^{-\alpha/2+1/(2p)-3/4}
\\[.25cm]
e_n^\deno(I,B_X)&\le& c2^{-\alpha k_0}\le cn^{-\alpha/2}.
\end{eqnarray*}
Combined with \eqref{AE0}--\eqref{AF9} and \eqref{J9} this completes the proof with $n_0=\max(2^{c(1)},c(3)2^{2c(2)})$.

\end{proof}
The following is a direct consequence of Proposition \ref{pro:6}.
\begin{corollary}
\label{cor:6} Let $1\le p<2<u\le \infty$ and $\alpha>1$. Then there are constants $c_1,c_2>0$ and $n_0\in\N$ such that for $n\ge n_0$
\begin{eqnarray*}
\lefteqn{c_1n^{1/(2u)-1/4}+c_1n^{1/4-1/(2p)}\le \frac{e_n^\ran(I,B_X)}{e_n^\ranno(I,B_X)}}
\notag\\[.2cm]
&\le& c_2n^{1/(2u)-1/4}(\log(n+1))^{1-1/u}+c_2n^{1/4-1/(2p)}(\log(n+1))^{1/2}.
\end{eqnarray*}
For $p=1$, $u=\infty$ this gives 
\begin{equation*}
c_1n^{1/4}(\log(n+1))^{-1}\le \frac{e_n^\ranno(I,B_X)}{e_n^\ran(I,B_X)}\le c_2n^{1/4}. 
\end{equation*}
\end{corollary}
\subsection{Infinite diagonal problems}
\label{sec:5.2}

A similar technique as used above can be applied to the results of \cite{Hei23c} on approximation in finite dimensional spaces. This analysis  leads firstly to an infinite dimensional problem with largest so far known gap of order $n^{1/2}(\log n)^{-1}$  for all $n$. Secondly, we can produce examples  with either the source space or the target space being Hilbertian and the gap being of order $n^{1/4}(\log n)^{-1/2}$. This shows that the adaption problem has a negative solution even in such situations. All this relates to the setting of standard information. For the case of linear information we refer to \cite{KNW23}, where an infinite dimensional problem with a Hilbertian target space and a gap of $n^{1/2}(\log n)^{-1/2}$ was found. 

For our purpose we have to consider a slightly different from \eqref{MA3} situation, namely that the problems 
$\mathcal{P}_k$ have different target spaces, so let $T_k\in L(X_k,G_k)$, where $X_k$ and $G_k$ are Banach spaces, let $\emptyset\ne \Lambda_k\subseteq X_k^*$, and put $\mathcal{P}_k=(B_{X_k},G_k,T_k,\K,\Lambda_k)$.
Let $1\le q_1\le \infty$ and define 
$$
G=\bigg(\bigoplus_{k=0}^\infty G_k\bigg)_{\ell_{q_1}}.
$$ 
We consider the diagonal problem $\mathcal{P}=(B_X,G,T,\K,\Lambda)$, where $X$ and $\Lambda$ are given as above in \eqref{CF7} and \eqref{CF3}, while  
\begin{eqnarray*}
T:X\to G,\quad T((x_k)_{k=0}^\infty)=(T_kx_k)_{k=0}^\infty.
\end{eqnarray*}
With 
\begin{equation}
\label{CF6}
\frac{1}{w}=\left(\frac{1}{q_1}-\frac{1}{p_1}\right)_+
\end{equation}
(thus $w=\infty$ for $p_1\le q_1$), the operator $T$ is well-defined and bounded iff
\begin{eqnarray}
\|(\|T_k\|)\|_{\ell_w}<\infty,
\label{CH7}
\end{eqnarray}
because
\begin{eqnarray*}
\|T\|&=&\sup_{\|(\|x_k\|)\|_{\ell_{p_1}}\le 1}\|(\|T_k(x_k)\|_{G_k})\big\|_{\ell_{q_1}}
=\sup_{\|(b_k)\|_{\ell_{p_1}}\le 1}\|(\sup_{\|x_k\|\le b_k}\|T_k(x_k)\|_{G_k})\big\|_{\ell_{q_1}}
\nonumber\\
&=&\sup_{\|(b_k)\|_{\ell_{p_1}}\le 1}\|(\|T_k\|b_k)\big\|_{\ell_{q_1}}=\|(\|T_k\|)\|_{\ell_w},
\end{eqnarray*}
where the latter equality is obvious for the case $p_1\le q_1$, while for $p_1>q_1$ we refer to \cite{Pie78}, the case $n=1$ of Th. 11.11.4. Analogously,  we have 
\begin{eqnarray}
\big\|T-T\sum_{k=0}^{k_1} P_k\big\|_{L(X,G)}&=& \big\|( \|T_k\|)_{k=k_1+1}^\infty\big\|_{\ell_w}\quad(k_1\in\N_0).
\label{AF1}
\end{eqnarray}
It turns out that we can use Lemma \ref{lem:6} also for the present situation to obtain
\begin{corollary}
\label{cor:2}
Let ${\rm set}\in\{\de,\deno,\ran,\ranno\}$ and suppose that \eqref{CH7} holds. Then for $n\in\N_0$
\begin{eqnarray}
\label{CF8}
 e_n^\set(T_k,B_{X_k}, G_k)\le e_n^\set(T,B_X,G).
\end{eqnarray}
Moreover, for $k_1\in\N_0$, $(n_k)_{k=0}^{k_1}\subset\N_0$, $n$ defined by $n=\sum_{k=0}^{k_1}n_k$, and $w$ given by \eqref{CF6},
\begin{eqnarray}
\label{CF9}
e_n^\set(T,B_X,G)\le \big\|( \|T_k\|)_{k=k_1+1}^\infty\big\|_{\ell_w}+\sum_{k=0}^{k_1} e_{n_k}^\set(T_k,B_{X_k}, G_k).
\end{eqnarray}
\end{corollary}
\begin{proof}
Define
\begin{eqnarray*}
&& U_k:G_k\to G, \quad U_kg_k=(\delta^{kl}g_k)_{l=0}^\infty,
\\
&&V_k:G\to G_k, \quad V_k(g_l)_{l=0}^\infty=g_k.
\end{eqnarray*}
We apply  Lemma \ref{lem:6} with $S_k=U_kT_k$ and $S=T$. We have 
$V_kS_k=T_k$,
consequently
\begin{eqnarray*}
e_n^\set(S_k,B_{X_k},G)&\le& \|U_k\|e_n^\set(T_k,B_{X_k},G_k)=e_n^\set(T_k,B_{X_k},G_k)
\\
&\le& \|V_k\| e_n^\set(S_k,B_{X_k},G)=e_n^\set(S_k,B_{X_k},G),
\end{eqnarray*}
thus 
\begin{equation}
e_n^\set(S_k,B_{X_k},G)=e_n^\set(T_k,B_{X_k},G_k).
\label{CG5}
\end{equation}
Now  \eqref{PA2} and \eqref{PA3} of Lemma \ref{lem:6} together with \eqref{AF1} and \eqref{CG5} give
\begin{eqnarray*}
 e_n^\set(T_k,B_{X_k}, G_k)=e_n^\set(S_k,B_{X_k},G)\le e_n^\set(S,B_X,G)= e_n^\set(T,B_X,G)
\end{eqnarray*}
and
\begin{eqnarray*}
e_n^\set(T,B_X,G)&=& e_n^\set(S,B_X,G)\le \big\|S-S\sum_{k=0}^{k_1} P_k\big\|_{L(X,G)}+\sum_{k=0}^{k_1} e_{n_k}^\set(S_k,B_{X_k}, G)
\\
&=& \big\|T-T\sum_{k=0}^{k_1} P_k\big\|_{L(X,G)}+\sum_{k=0}^{k_1} e_{n_k}^\set(T_k,B_{X_k}, G)
\\
&=& \big\|( \|T_k\|)_{k=k_1+1}^\infty\big\|_{\ell_w}+\sum_{k=0}^{k_1} e_{n_k}^\set(T_k,B_{X_k},G).
\end{eqnarray*}

\end{proof}

 Let $1\le p,q,u,v\le\infty$, $M_1,M_2\in\N$,  define 
\begin{equation*}
\label{AE7}
J^{M_1,M_2}:L_p^{M_1}\big(L_u^{M_2}\big)\to L_q^{M_1}\big(L_v^{M_2}\big), \quad J^{M_1,M_2}f=f,
\end{equation*}
and consider the problem 
$\big(B_{L_p^{M_1}(L_u^{M_2})},L_q^{M_1}(L_v^{M_2}),J^{M_1,M_2},\K,\Lambda^{M_1,M_2}\big)$, where $\Lambda^{M_1,M_2}$ is standard information as defined in \eqref{CI7}. Then we have the norm bound  
\begin{equation}
\label{AG5}
\|J^{M_1,M_2}:L_p^{M_1}\big(L_u^{M_2}\big)\to L_q^{M_1}\big(L_v^{M_2}\big)\|=M_1^{(1/p-1/q)_+}M_2^{(1/u-1/v)_+}
\end{equation}
and the following theorem, which is a part of Theorem 1 in \cite{Hei23c}. 
\begin{theorem}
\label{theo:2}
Let $1 \leq p,q,u,v \le \infty $, $p<q$, $u>v$.
Then there exist constants $0<c_0<1$, $c_{1-4}>0$,  such that for all $n,M_1,M_2\in{\mathbb{N}}$ with 
$n < c_0M_1M_2$ the following hold: 
\begin{eqnarray}
\label{AE8}
&&c_1M_1^{1/p-1/q}\Bigg(\left\lceil\frac{n}{M_1}\right\rceil^{1/u-1/v}+ \left\lceil\frac{n}{M_2}\right\rceil^{1/q-1/p}\Bigg)
\le e_n^{\rm ran }\Big(J^{M_1,M_2},B_{L_p^{M_1}(L_u^{M_2})},L_q^{M_1}(L_v^{M_2})\Big) 
\nonumber\\[.2cm]
&\le& c_2M_1^{1/p-1/q}\Bigg(\left\lceil\frac{n}{M_1\log(M_1+M_2)}\right\rceil^{1/u-1/v}+ \left\lceil\frac{n}{M_2\log(M_1+M_2)}\right\rceil^{1/q-1/p}\Bigg),
\end{eqnarray}
\begin{eqnarray}
c_3M_1^{1/p-1/q}
  &\le& e_n^{\rm ran-non}\Big(J^{M_1,M_2},B_{L_p^{M_1}(L_u^{M_2})},L_q^{M_1}(L_v^{M_2})\Big) 
\le c_4M_1^{1/p-1/q}.
\label{AE9}
\end{eqnarray}
\end{theorem}
\medskip
 We use Corollary \ref{cor:2} in the following situation. Let $\alpha\in\R$, $\alpha>2$ and define for $k\in\N_0$ 
\begin{equation*}
N_k=2^k,\quad X_k=L_p^{N_k}(L_u^{N_k}),\quad G_k=L_q^{N_k}(L_v^{N_k}),
\quad T_k=2^{-\alpha k}J^{N_k,N_k},
\end{equation*}
and let $\mathcal{P}_k=(B_{L_p^{N_k}(L_u^{N_k})},L_q^{N_k}(L_v^{N_k}),2^{-\alpha k}J^{N_k,N_k},\K,\Lambda^{N_k,N_k}) $. We only consider the case $p<q$, $u>v$. 
Furthermore, let $1\le p_1,q_1\le\infty$,
\begin{equation*}
X=\bigg(\bigoplus_{k=0}^\infty L_p^{N_k}(L_u^{N_k})\bigg)_{p_1}, \quad G=\bigg(\bigoplus_{k=0}^\infty L_q^{N_k}(L_v^{N_k})\bigg)_{\ell_{q_1}}, 
\end{equation*}
\begin{equation*}
J:X\to G,\quad J((x_k))=(2^{-\alpha k}J^{N_k,N_k}x_k)\in G.
\end{equation*}
By \eqref{AG5} we have $\|J^{N_k,N_k}:L_p^{N_k}(L_u^{N_k})\to L_q^{N_k}(L_v^{N_k})\|= N_k^{1/p-1/q}$ and with $w$ given by \eqref{CF6}, we conclude from \eqref{CH7} that the operator $J$ is bounded,
\begin{eqnarray*}
\|J\| &=&\|(\|2^{-\alpha k}J^{N_k,N_k}\|)_{k=0}^\infty\|_w =\|(2^{(-\alpha +1/p-1/q)k})_{k=0}^\infty\|_w <\infty,
\end{eqnarray*}
and similarly, from \eqref{AF1} that there is a constant $c>0$ such that for all $k_1\in\N_0$
\begin{eqnarray}
\label{PA7}
\|J-\sum_{k=0}^{k_1} JP_k\|&=&\|(\|2^{-\alpha k}J^{N_k,N_k}\|)_{k=k_1+1}^\infty\|_w
\notag\\
& =&\|(2^{(-\alpha +1/p-1/q)k})_{k=k_1+1}^\infty\|_w\le c 2^{(-\alpha +1/p-1/q)k_1}.
\end{eqnarray}
\begin{proposition}
\label{pro:7} Let $1\le p,q,u,v,p_1,q_1\le\infty$, $p<q$, $u>v$, and let $\alpha>2$. Then there are constants $c_{1-4}>0$ and $n_0\in\N$ such that for $n\ge n_0$ 
\begin{eqnarray*}
\lefteqn{c_1n^{(-\alpha+1/p-1/q+1/u-1/v)/2}
+c_1n^{-\alpha/2}\le e_n^\ran(J,B_X,G)}
\notag\\[.2cm]
&\le&
c_2n^{(-\alpha+1/p-1/q+1/u-1/v)/2}(\log(n+1))^{1/v-1/u}
+c_2n^{-\alpha/2}(\log(n+1))^{1/p-1/q},
\end{eqnarray*}
\vspace{-.4cm}
\begin{eqnarray*}
 c_3n^{(-\alpha+1/p-1/q)/2}&\le& e_n^\ranno(J,B_X,G)\le e_n^\deno(J,B_X,G)
\\[.1cm]
&\le& 2 e_n^\de(J,B_X,G) \le c_4n^{(-\alpha+1/p-1/q)/2}. 
\end{eqnarray*}
\end{proposition}
\begin{proof} To obtain the lower bounds  we start as in the proof of Proposition \ref{pro:6}, where now  $0<c(0)<1$ is the constant $c_0$ from Theorem \ref{theo:2}. We use the same definitions and relations \eqref{PB5}--\eqref{PB4}.
From \eqref{CF8}, \eqref{AE9}, and \eqref{PB4} we obtain 
\begin{eqnarray*}
e_{m_k}^\ran(J,B_X,G)&\ge&  2^{-\alpha k}e_{m_k}^\ran(J^{N_k,N_k},B_{L_p^{N_k}(L_u^{N_k})},L_q^{N_k}(L_v^{N_k}))
\nonumber\\
&\ge& c2^{-\alpha k} \big(N_k^{1/p-1/q+1/v-1/u}m_k^{1/u-1/v}+ N_k^{2/p-2/q}m_k^{1/q-1/p}\big)
\notag \\
&\ge& c2^{(-\alpha+1/p-1/q+1/u-1/v)k}+c2^{-\alpha k}
\\
e_{m_k}^\ranno(J,B_X,G)&\ge&  2^{-\alpha k}e_{m_k}^{\rm ran-non }(J^{N_k,N_k},B_{L_p^{N_k}(L_u^{N_k})},L_q^{N_k}(L_v^{N_k}))
\nonumber\\
&\ge& c2^{-\alpha k}N_k^{1/p-1/q} \ge c2^{(-\alpha+1/p-1/q)k}.
\end{eqnarray*}
Arguing analogously to \eqref{CK2}--\eqref{CJ9}, we obtain for 
$n\ge 2^{c(1)}$
\begin{eqnarray}
e_n^\ran(J,B_X,G)&\ge& cn^{(-\alpha+1/p-1/q+1/u-1/v)/2}
+cn^{-\alpha/2}\label{AG9}
\\
e_n^\ranno(J,B_X,G)&\ge& c 2^{(-\alpha+1/p-1/q)k}\ge cn^{(-\alpha+1/p-1/q)/2}.
\label{CI3}
\end{eqnarray}

Next we show the upper estimates. Also this part is largely analogous to
the respective part of the proof of Proposition \ref{pro:6}. So we use definitions and relations \eqref{AE3}--\eqref{PA6}, while \eqref{PC0} is replaced by
\begin{eqnarray}
\label{CH3}
\lefteqn{e_{n_k}^\ran(J^{N_k,N_k},B_{L_p^{N_k}(L_u^{N_k})}, L_q^{N_k}(L_v^{N_k}))}
\notag\\
&=&e_{n_k}^\de(J^{N_k,N_k},B_{L_p^{N_k}(L_u^{N_k})}, L_q^{N_k}(L_v^{N_k}))=0\qquad(0\le k< k_0).
\end{eqnarray}
Then we have by \eqref{PA0}--\eqref{CJ8}, \eqref{AG5},  and \eqref{AE8} for $k_0\le k\le k_1$
\begin{eqnarray}
\lefteqn{e_{n_k}^\ran(J^{N_k,N_k},B_{L_p^{N_k}(L_u^{N_k})}, L_q^{N_k}(L_v^{N_k}))}
\nonumber\\
&\le& c (k+1)^{1/v-1/u}N_k^{1/p-1/q+1/v-1/u}n_k^{1/u-1/v}+c (k+1)^{1/p-1/q}N_k^{2/p-2/q}n_k^{1/q-1/p}
\nonumber\\[.2cm]
&\le& c (k+1)^{1/v-1/u}2^{(1/p-1/q+1/v-1/u)k+(1/u-1/v)(2k_0-\delta(k-k_0))}
\nonumber\\  [.1cm]   
&&+c (k+1)^{1/p-1/q}2^{(2/p-2/q)k+(1/q-1/p)(2k_0-\delta(k-k_0))}      
 \label{CH4}
\\[.2cm]
\lefteqn{e_{n_k}^\de(J^{N_k,N_k},B_{L_p^{N_k}(L_u^{N_k})}, L_q^{N_k}(L_v^{N_k}))
\le N_k^{1/p-1/q}
=  2^{(1/p-1/q)k}.} \label{AG6}
\end{eqnarray}
We set
\begin{equation}
\label{CG9}
\beta=\frac{\alpha}{\alpha-1}, 
\end{equation}
so $\beta$ satisfies \eqref{AE3}, and the requirement \eqref{CI9} on $\delta>0$ turns into 
\begin{equation}
0<\delta<\alpha-2.\label{CK9}
\end{equation}
We get from \eqref{PA6},  \eqref{CF9}, \eqref{PA7}, and \eqref{CH3}--\eqref{CK9}
\begin{eqnarray}
\lefteqn{e_{c(3) 2^{2k_0}}^\ran(J,B_X,G)}
\nonumber\\
&\le& c2^{(-\alpha+1/p-1/q) k_1}
+c \sum_{k=k_0}^{k_1}2^{-\alpha k}(k+1)^{1/v-1/u}2^{(1/p-1/q+1/v-1/u)k+(1/u-1/v)(2k_0-\delta(k-k_0))}
\nonumber\\
&&+c \sum_{k=k_0}^{k_1}2^{-\alpha k}(k+1)^{1/p-1/q}2^{(2/p-2/q)k+(1/q-1/p)(2k_0-\delta(k-k_0))} 
\nonumber\\
&\le&  c2^{(-\alpha+1)\beta k_0}
\nonumber\\
&& +c (k_0+1)^{1/v-1/u}2^{(-\alpha+1/p-1/q+1/u-1/v)k_0}\sum_{k=k_0}^{k_1}2^{(-\alpha+1/p-1/q+(1/v-1/u)(1+\delta))(k-k_0)}
\nonumber\\
&&+c(k_0+1)^{1/p-1/q}2^{-\alpha k_0}\sum_{k=k_0}^{k_1}2^{(-\alpha+(1/p-1/q)(2+\delta))(k-k_0)}
\nonumber\\
&\le& c2^{-\alpha k_0}+c(k_0+1)^{1/v-1/u}2^{(-\alpha+1/p-1/q+1/u-1/v) k_0}+(k_0+1)^{1/p-1/q}2^{-\alpha k_0}
\nonumber\\
&\le& c_(k_0+1)^{1/v-1/u}2^{(-\alpha+1/p-1/q+1/u-1/v) k_0}+c(k_0+1)^{1/p-1/q}2^{-\alpha k_0}.
\label{CH6}
\end{eqnarray}
and 
\begin{eqnarray}
e_{c(3) 2^{2k_0}}^\de(J,B_X,G)
&\le& c2^{(-\alpha+1/p-1/q) k_1}
+c \sum_{k=k_0}^{k_1}2^{(-\alpha +1/p-1/q)k}\le c2^{(-\alpha+1/p-1/q)k_0}.
\label{AG7}
\end{eqnarray}
Now let $n\ge c(3)2^{2c(2)}$, then we can choose $k_0\ge c(1)$ so that 
\begin{equation*}
c(3)2^{2k_0}\le n<c(3)2^{2(k_0+1)}.
\end{equation*}
Together with \eqref{CH6} and \eqref{AG7} it follows that for $n\ge c(3)2^{2c(2)}$
\begin{eqnarray}
e_n^\ran(J,B_X,G)&\le&   cn^{(-\alpha+1/p-1/q+1/u-1/v)/2}(\log(n+1))^{1/v-1/u}\notag\\
&&+cn^{-\alpha/2}(\log(n+1))^{1/p-1/q}\qquad 
\\[.2cm]
e_n^\ranno(J,B_X,G)&\le& e_n^\deno(J,B_X,G) \le 2 e_n^\de(J,B_X,G)\le cn^{(-\alpha+1/p-1/q)/2},\label{CH8}
\end{eqnarray}
where we also used  \eqref{AF7} and \eqref{J9} in \eqref{CH8}. Combining this with \eqref{AG9} and \eqref{CI3}  concludes the proof.

\end{proof}
\begin{corollary}
\label{cor:5} Assume that $\alpha>2$. Then there are constants $c_{1-6}>0$ and $n_0\in\N$ such that in the case $p=1$, $u=\infty$, $q=\infty$, $v=1$
\begin{equation*}
c_1n^{1/2}(\log(n+1))^{-1}\le \frac{e_n^\ranno(J,B_X,G)}{e_n^\ran(J,B_X,G)}\le c_2n^{1/2}\quad (n\ge n_0). 
\end{equation*}
Furthermore, if $p=u=p_1=2$, $q=\infty$, $v=1$, then $X$ is a Hilbert space and
\begin{equation*}
c_3n^{1/4}(\log(n+1))^{-1/2}\le \frac{e_n^\ranno(J,B_X,G)}{e_n^\ran(J,B_X,G)}\le c_4n^{1/4}\quad (n\ge n_0).
\end{equation*}
 If $p=1$, $u=\infty$,  $q=v=q_1=2$, then $G$ is a Hilbert space and we have 
\begin{equation*}
c_5n^{1/4}(\log(n+1))^{-1/2}\le \frac{e_n^\ranno(J,B_X,G)}{e_n^\ran(J,B_X,G)}\le c_6n^{1/4}\quad (n\ge n_0). 
\end{equation*}
\end{corollary}
\begin{proof}
It follows from Proposition \ref{pro:7} that for $n\ge n_0$
\begin{eqnarray*}
&&c_1n^{(1/u-1/v)/2}
+c_1n^{(1/q-1/p)/2}\le \frac{e_n^\ran(J,B_X,G)}{e_n^\ranno(J,B_X,G)}
\\[.1cm]
&\le& c_2n^{(1/u-1/v)/2}(\log(n+1))^{1/v-1/u}
+c_2n^{(1/q-1/p)/2}(\log(n+1))^{1/p-1/q}.
\end{eqnarray*}

\end{proof}
\noindent {\bf Acknowledgement.} I thank Erich Novak and Robert Kunsch for various discussions on the topic of this paper, in particular during the  Dagstuhl Seminar 23351 ''Algorithms and Complexity for Continuous Problems''.

\end{document}